# PERTE D'INFORMATION DANS LES TRANSFORMATIONS DU JEU DE PILE OU FACE

By Jean Brossard et Christophe Leuridan

*Institut Fourier*

Soit $(\varepsilon_n)_{n \in \mathbf{Z}}$ un jeu de pile ou face, c'est-à-dire une suite de variables aléatoires indépendantes de loi $(\delta_{-1} + \delta_1)/2$, et $(H_n)_{n \in \mathbf{Z}}$ un processus à valeurs dans $\{-1, 1\}$, prévisible dans la filtration naturelle de $(\varepsilon_n)_{n \in \mathbf{Z}}$. Alors $(H_n \varepsilon_n)_{n \in \mathbf{Z}}$ est encore un jeu de pile ou face, dont la filtration naturelle est contenue dans celle de $(\varepsilon_n)_{n \in \mathbf{Z}}$. Le but de l'article est d'obtenir des conditions pour que ces filtrations soient égales et de décrire l'écart entre ces filtrations lorsqu'elles sont différentes. Nous nous intéressons plus particulièrement au cas des transformations homogènes, où le processus $(H_n \varepsilon_n)_{n \in \mathbf{Z}}$ est une fonctionnelle de $(\varepsilon_n)_{n \in \mathbf{Z}}$ qui commute avec les translations. Nous étudions de façon approfondie les transformations homogènes de longueur finie, où $H_n$ est de la forme $\phi(\varepsilon_{n-d}, \ldots, \varepsilon_{n-1})$ avec $d \in \mathbf{N}$ et $\phi \colon \{-1; 1\}^d \to \{-1; 1\}$ fixés.

**0. Introduction.** Soit $(\varepsilon_n)_{n \in \mathbf{Z}}$ un jeu de pile ou face, c'est-à-dire une suite de variables aléatoires indépendantes de loi $(\delta_{-1} + \delta_1)/2$, et $(H_n)_{n \in \mathbf{Z}}$ un processus à valeurs dans $\{-1, 1\}$, prévisible dans la filtration naturelle de $(\varepsilon_n)_{n \in \mathbf{Z}}$, notée $(\mathcal{F}_n^{\varepsilon})_{n \in \mathbf{Z}}$. On vérifie facilement que $(H_n \varepsilon_n)_{n \in \mathbf{Z}}$ est encore un jeu de pile ou face, dont la filtration naturelle est contenue dans celle de $(\varepsilon_n)_{n \in \mathbf{Z}}$. Plus précisément, pour tout $n \in \mathbf{Z}$,

- $H_n \varepsilon_n$ a pour loi $(\delta_{-1} + \delta_1)/2$;
- $H_n \varepsilon_n$ est indépendante de $\mathcal{F}_{n-1}^{\varepsilon}$;
- $\mathcal{F}_n^{\varepsilon} = \mathcal{F}_{n-1}^{\varepsilon} \vee \sigma(H_n \varepsilon_n)$.

Comme la tribu asymptotique $\mathcal{F}_{-\infty}^{\varepsilon}$ est triviale (d'après la loi du 0–1 de Kolmogorov), on pourrait croire comme l'ont fait Kallianpur et Wiener en 1956 (voir [3] pour une discussion de ce point) que l'égalité $\mathcal{F}_{\infty}^{\varepsilon} = \mathcal{F}_{\infty}^{H\varepsilon} \vee \mathcal{F}_n^{\varepsilon}$ entraîne forcément $\mathcal{F}_{\infty}^{\varepsilon} = \mathcal{F}_{\infty}^{H\varepsilon}$ en faisant tendre $n$ vers $-\infty$. En fait, il n'en n'est rien et la "propriété d'échange" (de l'intersection sur $n$ avec le suprémum des tribus) n'a pas lieu en général. La propriété d'échange a été étudiée de façon générale par H. von Weizsäcker [6].

Pour certaines transformations, la suite d'innovations $(H_n \varepsilon_n)_{n \in \mathbf{Z}}$ contient moins d'information que la suite $(\varepsilon_n)_{n \in \mathbf{Z}}$. Le contre-exemple le plus simple est celui où $H_n = \varepsilon_{n-1}$ pour tout $n \in \mathbf{Z}$. En effet, la suite $(\varepsilon_{n-1} \varepsilon_n)_{n \in \mathbf{Z}}$ est une fonctionnelle paire de la suite $(\varepsilon_n)_{n \in \mathbf{Z}}$. La transformation $\Phi_1 : (x_n)_{n \in \mathbf{Z}} \mapsto (x_{n-1} x_n)_{n \in \mathbf{Z}}$ de $\{-1; 1\}^{\mathbf{Z}}$ dans $\{-1; 1\}^{\mathbf{Z}}$ perd donc au moins un bit d'information puisque deux suites opposées ont même image par $\Phi_1$. En fait, on vérifie facilement que $\Phi_1$ perd exactement un bit d'information puisque pour tout $n \in \mathbf{Z}$, $\varepsilon_n$ est indépendante de $\sigma((\varepsilon_{k-1} \varepsilon_k)_{k \leq n})$ et $\sigma(\varepsilon_n) \vee \sigma((\varepsilon_{k-1} \varepsilon_k)_{k \leq n}) = \sigma((\varepsilon_k)_{k \leq n})$.

Dans cet article, nous nous intéressons aux problèmes suivants:

- A quelle condition la filtration $(\mathcal{F}_n^{H\varepsilon})_{n \in \mathbf{Z}}$ est-elle égale à la filtration $(\mathcal{F}_n^{\varepsilon})_{n \in \mathbf{Z}}$?
- Existe-t-il un complément indépendant, c'est-à-dire une variable aléatoire $U$ indépendante de $\mathcal{F}_\infty^{\varepsilon}$ telle que $\mathcal{F}_\infty^{H\varepsilon} \vee \sigma(U) = \mathcal{F}_\infty^{\varepsilon}$?
- L'information contenue dans $U$ est-elle équivalente à l'information contenue dans un nombre fini de bits indépendants?

Nous verrons dans la première partie que la réponse à la deuxième question est positive, et que la variable aléatoire $U$ est soit uniforme sur un ensemble fini, soit diffuse (et peut être alors choisie uniforme sur $[0, 1]$). En revanche, le nombre de valeurs de $U$ n'est pas toujours une puissance de 2 et il peut être infini, sauf dans le cas d'une transformation homogène de longueur finie, c'est-à-dire lorsque $H_n$ est de la forme $\phi(\varepsilon_{n-d}, \ldots, \varepsilon_{n-1})$ avec $d \in \mathbf{N}$ et $\phi : \{-1; 1\}^d \to \{-1; 1\}$ fixés.

Ce cas particulier est bien élucidé. Nous disposons d'une condition nécessaire et suffisante pour que la filtration $(\mathcal{F}_n^{H\varepsilon})_{n \in \mathbf{Z}}$ soit égale à la filtration $(\mathcal{F}_n^{\varepsilon})_{n \in \mathbf{Z}}$. Ce résultat recoupe en partie ceux de Rosenblatt [5] pour les chaînes de Markov "uniformes" sur un ensemble fini. Mais notre étude va plus loin puisque nous décrivons complètement l'information perdue lorsqu'il y a perte d'information, montrant en particulier que le nombre de valeurs que prend la variable aléatoire $U$ est toujours une puissance de 2. Autrement dit, on peut coder l'information perdue par un nombre fini de bits indépendants.

La condition de conservation de l'information que nous donnons pour une transformation homogène de longueur finie fournit un algorithme effectif. Dans l'article [5], Rosenblatt avait proposé une "procédure de calcul" pour déterminer s'il y a ou non perte d'information. Mais la méthode proposée est clairement très coûteuse en temps de calcul, et ne permettait guère d'explorer un grand nombre de cas. Nous remercions Roland Bacher d'avoir conçu un algorithme astucieux et de l'avoir fait tourner pour nous. Cet algorithme a permis de résoudre complètement le cas où $d = 4$ en donnant la liste complète des 782 applications $\phi : \{-1; 1\}^4 \to \{-1; 1\}$ (parmi $2^{16} = 65\ 536$) pour lesquelles la transformation associée perd de l'information. A partir de $d = 5$, il devient trop long de regarder les $2^{32} = 4\ 294\ 967\ 296$ applications $\phi \colon \{-1; 1\}^5 \to \{-1; 1\}$, mais il reste possible de déterminer rapidement si la transformation associée à une application $\phi$ fait perdre de l'information. Ces résultats obtenus par l'informatique nous ont été précieux pour trouver des exemples intéressants et guider nos recherches.

**1. Résultats généraux sur l'information perdue.** Toutes les tribus que nous considérons dans cet article sont supposées complètes. Nous utiliserons la notation suivante: pour $x = (x_n)_{n \in \mathbf{Z}} \in \{-1; 1\}^{\mathbf{Z}}$, et $n \in \mathbf{Z}$, on note $x_{n]}$ (qu'on lit "$x$ avant $n$") la suite $(x_{n+k})_{k \in -\mathbf{N}}$.

Le but de cette partie est de décrire l'écart entre les filtrations $(\varepsilon_n)_{n \in \mathbf{Z}}$ et $(H_n \varepsilon_n)_{n \in \mathbf{Z}}$ lorsque $(H_n)_{n \in \mathbf{Z}}$ est un processus prévisible dans la filtration naturelle de $(\varepsilon_n)_{n \in \mathbf{Z}}$.

*1.1. Le résultat principal.* Soit $(H_n)_{n \in \mathbf{Z}}$ un processus prévisible dans la filtration naturelle de $(\varepsilon_n)_{n \in \mathbf{Z}}$. Pour tout $n \in \mathbf{Z}$, $H_n = h_n(\varepsilon_{n-1]})$ où $h_n$ est une application mesurable de $\{-1; 1\}^{-\mathbf{N}}$ dans $\{-1; 1\}$.

Notons $\Phi \colon \{-1; 1\}^{\mathbf{Z}} \to \{-1; 1\}^{\mathbf{Z}}$ la transformation définie par $\Phi(x) = (h_n(x_{n-1]})x_n)_{n \in \mathbf{Z}}$ pour tout $x = (x_n)_{n \in \mathbf{Z}} \in \{-1; 1\}^{\mathbf{Z}}$. Autrement dit, $\Phi$ est l'application qui fait passer du jeu de pile ou face $(\varepsilon_n)_{n \in \mathbf{Z}}$ au jeu de pile ou face $(H_n \varepsilon_n)_{n \in \mathbf{Z}}$.

Le résultat principal que nous allons démontrer est le suivant:

THÉORÈME 1. *La loi de $\varepsilon$ sachant $\Phi(\varepsilon)$ est soit presque sûrement diffuse, soit presque sûrement uniforme sur un ensemble fini de cardinal constant.*

La démonstration du théorème repose sur les deux faits suivants:

LEMME 2. *Soit $(X_n)_{n \in \mathbf{N}^*}$ une suite de variables aléatoires indépendantes de même loi $\mu$. Alors la suite des probabilités $P[X_1 = \cdots = X_d]$ détermine les masses des atomes de $\mu$.*

DÉMONSTRATION. En effet, soient $p_1 \geq p_2 \geq \cdots$ les masses des atomes de $\mu$ rangées dans l'ordre décroissant (avec $p_n = 0$ si $\mu$ a moins de $n$ atomes) et $q_d = P[X_1 = \cdots = X_d]$. Alors $q_d = \sum_n p_n^d$ pour tout $d \geq 2$. Quand $d \to +\infty$, on a donc $q_d^{1/d} \to p_1$, $(q_d - p_1^d)^{1/d} \to p_2$. $\square$

LEMME 3. *Si deux suites $x, y \in \{-1; 1\}^{\mathbf{Z}}$ vérifient $\Phi(x) = \Phi(y)$ et coïncident jusqu'à un certain rang, elles sont égales. Autrement dit, pour tout $z \in \{-1; 1\}^{\mathbf{Z}}$, la restriction à $\Phi^{-1}(\{z\})$ de l'application $x \mapsto x_{n]}$ est injective.*

DÉMONSTRATION. En effet, si $\Phi(x) = \Phi(y) = z$ et si $x_{n]} = y_{n]}$, alors $x_{n+1} = y_{n+1}$ puisque $h_{n+1}(x_{n]})x_{n+1} = z_{n+1} = h_{n+1}(y_{n]})y_{n+1}$. On en déduit le résultat par récurrence. □

Passons maintenant à la démonstration du théorème.

DÉMONSTRATION DU THÉORÈME 1. Soit $\mu = \bigotimes_{n \in \mathbf{Z}} \frac{\delta_{-1} + \delta_1}{2}$ la loi uniforme $\{-1; 1\}^{\mathbf{Z}}$. Pour tout $\alpha \in \{-1; 1\}^{\mathbf{Z}}$, on note $\nu(\alpha, \cdot)$ la loi de $\varepsilon$ sachant que $\Phi(\varepsilon) = \alpha$. Pour tout $n \in \mathbf{Z}$, on note $\nu_n(\alpha, \cdot)$ l'image de $\nu(\alpha, \cdot)$ par l'application $x \mapsto x_{n]}$, qui est donc la loi de $\varepsilon_{n]}$ sachant que $\Phi(\varepsilon) = \alpha$. Enfin pour tout $\beta \in \{-1; 1\}^{-\mathbf{N}}$, notons $\nu'_n(\beta, \cdot)$ la loi de $\varepsilon_{n]}$ sachant que $\Phi(\varepsilon)_{n]} = \beta$.

La démonstration du théorème s'effectue en trois étapes.

*Étape* 1. Pour $\mu$-presque tout $\alpha \in \{-1; 1\}^{\mathbf{Z}}$ et pour tout $n \in \mathbf{Z}$, les atomes de $\nu_n(\alpha, \cdot)$ sont les images par l'application $x \mapsto x_{n]}$ des atomes de $\nu(\alpha, \cdot)$, et il y a conservation des masses. De plus, $\nu_n(\alpha, \cdot) = \nu'_n(\alpha_{n]}, \cdot)$.

En effet, pour $\mu$-presque tout $\alpha \in \{-1; 1\}^{\mathbf{Z}}$, la loi $\nu(\alpha, \cdot)$ est portée par $\Phi^{-1}(\{\alpha\})$. Comme la restriction à $\Phi^{-1}(\{\alpha\})$ de l'application $x \mapsto x_{n]}$ est injective, d'après le lemme 3, les atomes de $\nu_n(\alpha, \cdot)$ sont les images des atomes de $\nu(\alpha, \cdot)$, et il y a conservation des masses.

La dernière égalité provient du fait que la loi de $\varepsilon_{n]}$ sachant $H\varepsilon$ est aussi la loi de $\varepsilon_{n]}$ sachant $(H\varepsilon)_{n]}$ par indépendance de $(H_k\varepsilon_k)_{k \geq n+1}$ et de $\mathcal{F}_n^\varepsilon$.

*Étape* 2. Les masses des atomes de $\nu(\alpha, \cdot)$ sont presque sûrement constantes.

En effet, soit $(\xi^n)_{n \in \mathbf{N}}$ une suite de jeux de pile ou face telle que conditionnellement à $\xi^0 = \alpha$, les variables aléatoires $(\xi^n)_{n \in \mathbf{N}^*}$ soient indépendantes de loi $\nu(\alpha, \cdot)$. Soit $d \geq 2$. Alors presque sûrement,

$$P[\xi^1 = \cdots = \xi^d | \xi^0] = P[\xi_{n]}^1 = \cdots = \xi_{n]}^d | \xi^0]$$

$$= P[\xi_{n]}^1 = \cdots = \xi_{n]}^d | \xi_{n]}^0],$$

la première égalité provenant du Lemme 3 et la seconde de l'égalité $\nu_n(\alpha, \cdot) = \nu'_n(\alpha_{n]}, \cdot)$ pour $\mu$-presque tout $\alpha$. Par conséquent, la variable aléatoire $P[\xi^1 = \cdots = \xi^d | \xi^0]$ est mesurable pour la tribu asymptotique $\mathcal{F}_{-\infty}^{\xi^0}$, donc presque sûrement constante. On conclut grâce au Lemme 2.

*Étape* 3. Les masses des atomes de $\nu(\alpha, \cdot)$ sont presque sûrement égales.

En effet, notons $p_1(\alpha)$ la plus grande des masses des atomes de $\nu(\alpha, \cdot)$. On sait d'après la deuxième étape que $p_1(\alpha)$ est constant pour presque tout $\alpha$. Notons $p_1$ cette constante. Notons $A(\alpha)$ l'ensemble des atomes de $\nu(\alpha, \cdot)$ de masse $p_1$ et $A_n(\alpha)$ l'ensemble des atomes de $\nu_n(\alpha, \cdot)$ de masse $p_1$. On sait d'après la première étape que pour $\mu$-presque tout $\alpha \in \{-1; 1\}^{\mathbf{Z}}$, la loi $\nu(\alpha, \cdot)$ est portée par $\Phi^{-1}(\{\alpha\})$ et les éléments de $A_n(\alpha)$ sont les images des éléments de $A(\alpha)$ par l'application $x \mapsto x_{n]}$. De plus, d'après le Lemme 3, la restriction à $\Phi^{-1}(\{\alpha\})$ de l'application $x \mapsto x_{n]}$ est injective. Donc pour tout $n \in \mathbf{Z}$,

$$[\varepsilon \in A(\Phi(\varepsilon))] = [\varepsilon_{n]} \in A_n(\Phi(\varepsilon))] \qquad \text{presque sûrement.}$$

Mais d'après la première étape, $A_n(\Phi(\varepsilon)) = A'_n(\Phi(\varepsilon)_{n]})$ presque sûrement, en notant $A'_n(\beta)$ l'ensemble des atomes de $\nu'_n(\beta, \cdot)$ de masse $p_1$. Donc l'événement $[\varepsilon_{n]} \in A_n(\Phi(\varepsilon))]$ appartient à $\mathcal{F}_n^\varepsilon$, ce qui montre que l'événement $[\varepsilon \in A(\Phi(\varepsilon))]$ appartient à $\mathcal{F}_{-\infty}^\varepsilon$, ce qui entraîne que $P[\varepsilon \in A(\Phi(\varepsilon))]$ vaut 0 ou 1. Par ailleurs, pour presque tout $\alpha \in \{-1; 1\}^{\mathbf{Z}}$,

$$P[\varepsilon \in A(\Phi(\varepsilon))|\Phi(\varepsilon) = \alpha] = \nu(\alpha, A(\alpha))$$
$$= p_1 \operatorname{Card}(A(\alpha)).$$

On a donc deux possibilités:

- Soit pour presque tout $\alpha$, $p_1 \operatorname{Card}(A(\alpha)) = 1$: la loi $\nu(\alpha, \cdot)$ est alors uniforme sur un ensemble de cardinal $1/p_1$.
- Soit pour presque tout $\alpha$, $p_1 \operatorname{Card}(A(\alpha)) = 0$: la loi $\nu(\alpha, \cdot)$ est alors diffuse. Le théorème est donc démontré.   $\square$

La démonstration du théorème 1 montre plus généralement que pour tout $n \in \mathbf{Z}$, la loi de $\varepsilon_{n]}$ sachant $\Phi(\varepsilon)_{n]}$ est soit presque sûrement diffuse, soit presque sûrement uniforme sur un ensemble fini de cardinal constant. En utilisant une injection borélienne de $\{-1; 1\}^{-\mathbf{N}}$ dans $[0, 1]$, des fonctions quantiles convenables et le fait que les atomes d'une probabilité sur $[0, 1]$ dépendent de façon mesurable de la probabilité, on en déduit le résultat suivant.

COROLLAIRE 4. *Pour tout $n \in \mathbf{Z}$, il existe une variable aléatoire $U_n$ mesurable pour $\mathcal{F}_n^\varepsilon$, indépendante de $\mathcal{F}_\infty^{H\varepsilon}$ et uniforme sur un ensemble fini (de cardinal indépendant de $n$) ou sur $[0, 1]$ telle que $\mathcal{F}_n^{H\varepsilon} \vee \sigma(U_n) = \mathcal{F}_n^\varepsilon$ et $\mathcal{F}_\infty^{H\varepsilon} \vee \sigma(U_n) = \mathcal{F}_\infty^\varepsilon$.*

REMARQUE. Une conséquence que nous utiliserons souvent est que l'information perdue s'ajoute par composition des transformations: si $\Phi$ fait perdre une information à $p$ valeurs et $\Psi$ une information à $q$ valeurs, alors la composée $\Psi \circ \Phi$ fait perdre une information à $pq$ valeurs.

Pour tout $d \in \mathbf{N}^*$ fixé, il est facile de construire une transformation faisant perdre $d$ bits d'information. Il suffit par exemple de considérer $\Phi : (x_n)_{n \in \mathbf{Z}} \mapsto (x_{n-d}x_n)_{n \in \mathbf{Z}}$. En effet, on vérifie facilement que pour tout $n \in \mathbf{Z}$, le $d$-uplet $(\varepsilon_{n-d+1}, \ldots, \varepsilon_n)$ est indépendant de $\sigma((\varepsilon_{k-d}\varepsilon_k)_{k \leq n})$ et $\sigma(\varepsilon_{n-d+1}, \ldots, \varepsilon_n) \vee \sigma((\varepsilon_{k-d}\varepsilon_k)_{k \leq n}) = \sigma((\varepsilon_k)_{k \leq n})$.

On peut se demander s'il est possible de perdre une information ayant un nombre fini de valeurs qui ne soit pas une puissance de 2 ou de perdre une information infinie. Nous allons voir que la réponse est positive, même dans le cas de transformations homogènes. Mais nous verrons dans la quatrième partie que la réponse est négative pour les transformations homogènes de longueur finie.

1.2. *Exemple de perte d'information ternaire.* On prend $\Omega = \mathbf{R}$, muni de la tribu $\mathcal{A}$ des boréliens 1-périodiques et de la probabilité $P$ qui à $A \in \mathcal{A}$ associe la mesure de Lebesgue de $A \cap [0,1[$. A chaque $\omega \in \Omega$, on associe le développement dyadique de sa partie fractionnaire:

$$\mathrm{Frac}(\omega) = \sum_{n=-\infty}^{0} 2^{n-1} D_n(\omega)$$

avec $D_n(\omega) \in \{0,1\}$ pour tout $n \in -\mathbf{N}$. Pour tout $n \in -\mathbf{N}$, on note $\varepsilon_n(\omega) = 2 D_n(\omega) - 1$. La suite $(\varepsilon_n)_{n \in -\mathbf{N}}$ est alors un jeu de pile ou face. Signalons que le fait que la suite soit indexée par $-\mathbf{N}$ au lieu de $\mathbf{Z}$ n'a pas d'importance puisque le problème se situe en $-\infty$.

Pour tout $n \in -\mathbf{N}$, la tribu $\mathcal{F}_n^\varepsilon = \sigma((\varepsilon_k)_{k \leq n})$ est la tribu des boréliens $2^n$-périodiques. Notons $\mathcal{G}_n$ la tribu des boréliens $\frac{2^n}{3}$-périodiques et $U$ la variable aléatoire $U$ valant 0 sur $[0, \frac{1}{3}[$, 1 sur $[\frac{1}{3}, \frac{2}{3}[$ et 2 sur $[\frac{2}{3}, 1[$. Alors $U$ est indépendante de $\mathcal{G}_0$, de loi uniforme sur $\{0, 1, 2\}$ et vérifie $\mathcal{G}_0 \vee \sigma(U) = \mathcal{F}_0$.

Par ailleurs, $(\mathcal{G}_n)_{n \leq 0}$ est la filtration naturelle de la suite $(\eta_n)_{n \leq 0}$ définie par $\eta_n(\omega) = \varepsilon_n(3\omega)$ pour tout $n \in -\mathbf{N}$ et pour tout $\omega \in \mathbf{R}$.

Ecrivons $\eta_n = H_n \varepsilon_n$ en posant $H_n(\omega) = \varepsilon_n(\omega)\varepsilon_n(3\omega)$ pour tout $\omega \in \mathbf{R}$. On voit que la variable aléatoire $H_n$ est mesurable pour $\mathcal{F}_{n-1}^\varepsilon$ en remarquant que pour tout $\omega \in \mathbf{R}$,

$$H_n(\omega + 2^{n-1}) = \varepsilon_n(\omega + 2^{n-1})\varepsilon_n(3\omega + 3 \cdot 2^{n-1})$$
$$= (-\varepsilon_n(\omega))(-\varepsilon_n(3\omega))$$
$$= H_n(\omega).$$

Et comme

$$H_n(\omega) = \varepsilon_0(2^{-n}\omega)\varepsilon_0(3 \cdot 2^{-n}\omega)$$
$$= H_0(2^{-n}\omega),$$

la transformation qui fait passer de la suite $(\varepsilon_n)_{n \in -\mathbf{N}}$ à la suite $(\eta_n)_{n \in -\mathbf{N}}$ est homogène. Plus précisément, on vérifie immédiatement que pour tout $\omega \in [0,1[$,

$$H_0(\omega) = -1 \quad \Longleftrightarrow \quad \omega \in [\tfrac{1}{6}, \tfrac{2}{6}[ \cup [\tfrac{4}{6}, \tfrac{5}{6}[$$
$$\Longleftrightarrow \quad \mathrm{Frac}(2\omega) \in [\tfrac{1}{3}, \tfrac{2}{3}[$$
$$\Longleftrightarrow \quad \sum_{n=-\infty}^{-1} 2^{n-1} D_{n-1}(\omega) \in [\tfrac{1}{3}, \tfrac{2}{3}[.$$

Par conséquent, $H_0(\omega) = h(\varepsilon_{-1]}(\omega)) = h((\varepsilon_{n-1})_{n \in -\mathbf{N}})$, en notant

$$h((x_n)_{n \in \mathbf{N}}) = \begin{vmatrix} -1, & \text{si } \sum_{n=-\infty}^{-1} 2^{n-1} \frac{1+x_n}{2} \in \left[\frac{1}{3}, \frac{2}{3}\right[, \\ 1, & \text{sinon.} \end{vmatrix}$$

Par périodicité cette dernière égalité reste valable pour tout $\omega \in \mathbf{R}$. Plus généralement on montre que $H_n(\omega) = H_0(2^{-n}\omega)$ pour tout $n \in -\mathbf{N}$ grâce à l'égalité $H_n(\omega) = h(\varepsilon_{n-1]}(\omega))$.

REMARQUE. Dans ce qui précède, le nombre 3 ne joue aucun rôle particulier et peut être remplacé par n'importe quel nombre impair. Par composition avec une transformation faisant perdre un nombre fini de bits, on obtient donc un exemple de transformation où l'information perdue est uniforme sur un ensemble fini de cardinal quelconque.

1.3. *Exemple de perte d'information infinie.* Un exemple simple de transformation faisant perdre une information infinie est la transformation $\Phi$ : $(x_n)_{n \in \mathbf{Z}} \mapsto (y_n)_{n \in \mathbf{Z}}$, avec $y_n = x_{2n} x_n$ pour $n \leq -1$ et $y_n = x_n$ pour $n \geq 0$. En effet, pour tout entier impair négatif $m$, la transformation $\Phi$ agit sur la sous-suite $(\ldots, x_{4m}, x_{2m}, x_m)$ comme la transformation $\Phi_1 : (x_n)_{n \in \mathbf{Z}} \mapsto (x_{n-1} x_n)_{n \in \mathbf{Z}}$. On perd donc un bit d'information pour chacune de ces sous-suites, qui sont indépendantes et en nombre infini. Mais cette transformation n'est pas homogène. Voyons maintenant un exemple de transformation homogène faisant perdre une information infinie.

Pour tout $k \in \mathbf{N}^*$ et pour toute suite $x \in \{-1; 1\}^{-\mathbf{N}}$, notons $h_k(x)$ le signe de la dernière répétition de longueur $k$ dans la suite $x$. Autrement dit, $h_k(x) = x_0$ si $x_0 = \cdots = x_{-k+1}$, $h_k(x) = x_{-1}$ si $x_0 \neq x_{-1} = \cdots = x_{-k}$, etc. . .

Soit $\Phi_k$ la transformation définie par $\Phi_k(x) = (h_k(x_{n-1]}) x_n)_{n \in \mathbf{Z}}$ pour tout $x \in \{-1; 1\}^{\mathbf{Z}}$. Notons $\varepsilon^k$ l'image de la suite $\varepsilon = (\varepsilon_n)_{n \in \mathbf{Z}}$ par la composée $\Phi_k \circ \cdots \circ \Phi_1$. Alors pour tout $n \in \mathbf{Z}$,

$$\varepsilon_n^k \varepsilon_{n-1}^k = \varepsilon_n^{k-1} \varepsilon_{n-1}^{k-1} h_k(\varepsilon_{n-1]}^{k-1}) h_k(\varepsilon_{n-2]}^{k-1}).$$

Or si $x \in \{-1; 1\}^{-\mathbf{N}}$, il faut que $x_0 = x_{-1} = \cdots = x_{-k+1} \neq x_k$ pour que $h_k(x)$ diffère de $h_k(x_{-1]})$. Donc

$$P[h_k(\varepsilon_{n-1]}^{k-1}) \neq h_k(\varepsilon_{n-2]}^{k-1})] \leq \frac{1}{2^k}.$$

On a donc presque sûrement $h_k(\varepsilon_{n-1]}^{k-1}) h_k(\varepsilon_{n-2]}^{k-1}) = 1$ pour tout $k$ assez grand, si bien que le produit $\varepsilon_n^k \varepsilon_{n-1}^k$ tend vers

$$\eta_n = \varepsilon_n \varepsilon_{n-1} \prod_{l=1}^{+\infty} h_l(\varepsilon_{n-1]}^{l-1}) h_l(\varepsilon_{n-2]}^{l-1}).$$

La suite $(\eta_n)_{n\in\mathbf{Z}}$ est bien l'image de la suite $(\varepsilon_n)_{n\in\mathbf{Z}}$ par une transformation homogène et cette transformation fait perdre une infinité de bits indépandants, puisque pour tout $k \in \mathbf{N}^*$, la suite $(\eta_n)_{n\in\mathbf{Z}}$ est une fonction de la suite $\varepsilon^k$ et la suite $\varepsilon^k$ est une fonction paire de la suite $\varepsilon^{k-1}$.

1.4. *Une condition nécessaire et suffisante de conservation de l'information.* Soit $\Phi$ la transformation définie par $\Phi(x) = (h_n(x_{n-1]})x_n)n \in \mathbf{Z}$ pour tout $x \in \{-1, 1\}^{\mathbf{Z}}$, où pour tout $n \in \mathbf{Z}$, $h_n$ est une application mesurable de $\{-1; 1\}^{-\mathbf{N}}$ dans $\{-1; 1\}$. L'image par $\Phi$ de la suite $(\varepsilon_n)_{n\in\mathbf{Z}}$ est donc la suite $(H_n\varepsilon_n)_{n\in\mathbf{Z}}$ où $H_n = h_n(\varepsilon_{n-1]})$ pour tout $n \in \mathbf{Z}$.

Notons $E = \{-1, 1\}^{-\mathbf{N}}$. Pour tout $n \in \mathbf{Z}$ et pour tout $\alpha \in \{-1, 1\} = \{-, +\}$, notons $f_\alpha^{(n)}$ l'application de $E$ dans lui-même définie par

$$f_\alpha^{(n)}((x_k)_{k\in-\mathbf{N}}) = (y_k)_{n\in-\mathbf{N}} \qquad \text{avec} \ \left| \begin{array}{ll} y_k = x_{k+1}, & \text{si } n \leq -1, \\ y_0 = \alpha h_n((x_k)_{k\in-\mathbf{N}}). \end{array} \right.$$

Par construction, les applications $f_+^{(n)}$ et $f_-^{(n)}$ vérifient

$$\varepsilon_{n]} = f_{H_n\varepsilon_n}^{(n)}(\varepsilon_{n-1]})$$

puisque $\varepsilon_n = H_n\varepsilon_n h_n(\varepsilon_{n-1]})$. Considérons les noyaux de transition sur $E^2$ définis par

$$N_n((a, b), \cdot) = \tfrac{1}{2}(\delta_{(f_+^{(n)}(a), f_+^{(n)}(b))} + \delta_{(f_-^{(n)}(a), f_-^{(n)}(b))}).$$

La forme des noyaux de transition entraîne la propriété suivante:

PROPOSITION 5. *Il y a équivalence entre:*

(1) $(X_n, Y_n)_{n\in\mathbf{Z}}$ *est une chaîne de Markov sur* $E^2$ *pour les noyaux de transition* $(N_n)_{n\in\mathbf{Z}}$;

(2) *il existe deux jeux de pile ou face* $(\xi_n)_{n\in\mathbf{Z}}$ *et* $(\eta_n)_{n\in\mathbf{Z}}$ *ayant même image* $(\zeta_n)_{n\in\mathbf{Z}}$ *par* $\Phi$ *tels que pour tout* $n \in \mathbf{Z}$, $\zeta_n$ *soit indépendant de* $(\xi_{n-1]}, \eta_{n-1]})$ *et* $(X_n, Y_n) = (\xi_{n]}, \eta_{n]})$ *presque sûrement.*

DÉMONSTRATION. L'implication $(2) \Rightarrow (1)$ est simple: il suffit de remarquer que si $(\zeta_n)_{n\in\mathbf{Z}}$ est l'image commune par $\Phi$ des jeux de pile ou face $(\xi_n)_{n\in\mathbf{Z}}$ et $(\eta_n)_{n\in\mathbf{Z}}$ alors $\xi_{n]} = f_{\zeta_n}^{(n)}(\xi_{n-1]})$ et $\eta_{n]} = f_{\zeta_n}^{(n)}(\eta_{n-1]})$.

La réciproque est à peine plus difficile: l'hypothèse (1) et la forme des noyaux de transition montrent que pour tout $n \in \mathbf{Z}$, il existe presque sûrement un unique $\zeta_n \in \{-1, 1\}$, tel que $X_n = f_{\zeta_n}^{(n)}(X_{n-1})$ et $Y_n = f_{\zeta_n}^{(n)}(Y_{n-1})$ et que la loi de $\zeta_n$ sachant $\mathcal{F}_{n-1}^{(X,Y)} = \sigma((X_k, Y_k)_{k\leq n-1})$ est la loi uniforme sur $\{-1, 1\}$. Par conséquent, $(\zeta_n)_{n\in\mathbf{Z}}$ est un jeu de pile ou face dans la filtration $(\mathcal{F}_n^{(X,Y)})_{n\in\mathbf{Z}}$. Notons $\xi_n$ et $\eta_n$ les derniers termes des suites $X_n$ et $Y_n$. La forme des applications $f_+^{(n)}$ et $f_-^{(n)}$ entraîne que presque sûrement, pour tout $n \in \mathbf{Z}$, $X_n = \xi_{n]}$, $Y_n = \eta_{n]}$. Les égalités $\xi_n = h_n(\xi_{n-1]})\zeta_n$ et $\eta_n = h_n(\eta_{n-1]})\zeta_n$ montrent que $(\xi_n)_{n \in \mathbf{Z}}$ et $(\eta_n)_{n \in \mathbf{Z}}$ sont aussi des jeux de pile ou face dans la filtration $(\mathcal{F}_n^{(X,Y)})_{n \in \mathbf{Z}}$, dont l'image par $\Phi$ est la suite $(\zeta_n)_{n \in \mathbf{Z}}$.  $\square$

Nous pouvons donner une condition nécessaire et suffisante de conservation de l'information. Cette condition, qui fait intervenir un couplage, ressemble à la condition de confort introduite par Emery et Schachermayer [1].

Théorème 6.  *Il y a équivalence entre:*

(1)  $\mathcal{F}_\infty^{H\varepsilon} = \mathcal{F}_\infty^\varepsilon$;
(2)  *pour tout* $n \in \mathbf{Z}$, $\mathcal{F}_n^{H\varepsilon} = \mathcal{F}_n^\varepsilon$;
(3)  *Les seules chaînes de Markov indexées par* $\mathbf{Z}$ *pour les noyaux de transition* $(N_n)_{n \in \mathbf{Z}}$ *vivent sur la diagonale de* $E^2$ *(presque sûrement).*

Remarque.  Notons $D$ la diagonale de $E^2$. Bien que les noyaux de transition $(N_n)_{n \in \mathbf{Z}}$ dépendent du temps, $D$ est une partie absorbante et les noyaux de transition induits sur $D$ ne dépendent pas du temps. La Proposition 5 montre que la seule chaîne portée par $D$ pour les noyaux de transition $(N_n)_{n \in \mathbf{Z}}$ est la chaîne stationnaire sur $D$.

Démonstration.  L'implication $(2) \Rightarrow (1)$ est triviale.

L'implication $(1) \Rightarrow (2)$ découle de l'inclusion $\mathcal{F}_n^{H\varepsilon} \subset \mathcal{F}_n^\varepsilon$ et de l'indépendance de $\mathcal{F}_n^\varepsilon$ et de $(H_k \varepsilon_k)_{k \geq n+1}$. Si $\mathcal{F}_\infty^{H\varepsilon} = \mathcal{F}_\infty^\varepsilon$, on peut en effet écrire, pour tout $A \in \mathcal{F}_n^\varepsilon$,

$$P[A|\mathcal{F}_n^{H\varepsilon}] = P[A|\mathcal{F}_n^{H\varepsilon} \vee \sigma((H_k \varepsilon_k)_{k \geq n+1})]$$
$$= P[A|\mathcal{F}_\infty^{H\varepsilon}] = P[A|\mathcal{F}_\infty^\varepsilon] = \mathbf{I}_A,$$

d'où $A \in \mathcal{F}_n^{H\varepsilon}$.

Montrons l'implication $(1) \Rightarrow (3)$. Si $\mathcal{F}_\infty^{H\varepsilon} = \mathcal{F}_\infty^\varepsilon$, il existe une application mesurable $F$ de $\{-1; 1\}^{\mathbf{Z}}$ dans $\{-1; 1\}^{\mathbf{Z}}$ telle que $\varepsilon = F(H\varepsilon) = F(\Phi(\varepsilon))$ presque sûrement, c'est-à-dire $F \circ \Phi = \mathrm{Id}$ $\mu$-presque partout. Soit alors $(X_n, Y_n)_{n \in \mathbf{Z}}$ une chaîne de Markov sur $E^2$ pour les noyaux de transition $(N_n)_{n \in \mathbf{Z}}$. D'après la Proposition 5, il existe deux jeux de pile ou face $(\xi_n)_{n \in \mathbf{Z}}$ et $(\eta_n)_{n \in \mathbf{Z}}$ ayant même image par $\Phi$ tels que presque sûrement, pour tout $n \in \mathbf{Z}$, $X_n = \xi_{n]}$ et $Y_n = \eta_{n]}$. On en déduit que les suites $(\xi_n)_{n \in \mathbf{Z}}$ et $(\eta_n)_{n \in \mathbf{Z}}$ sont égales presque sûrement en composant avec $F$. Par conséquent, les suites $(X_n)_{n \in \mathbf{Z}}$ et $(Y_n)_{n \in \mathbf{Z}}$ sont égales presque sûrement, ce qui montre que la chaîne $(X_n, Y_n)_{n \in \mathbf{Z}}$ est à valeurs dans $D$.

Montrons l'implication $(3) \Rightarrow (1)$. Notons à nouveau $\nu(\alpha, \cdot)$ et $\nu_n(\alpha, \cdot)$ la loi de $\varepsilon$ et la loi de $\varepsilon_{n]}$ sachant que $\Phi(\varepsilon) = \alpha$, et prenons trois jeux de pile ou face $\xi, \eta, \zeta$ tels que conditionnellement à $\zeta = \alpha$, les variables aléatoires $\xi$ et $\eta$ soient indépendantes de loi $\nu(\alpha, \cdot)$. Alors $\Phi(\xi) = \Phi(\eta) = \zeta$ presque sûrement.

De plus, pour tout $n \in \mathbf{Z}$, la variable aléatoire $\zeta_{n+1}$ est indépendante de $(\xi_{n]}, \eta_{n]})$. En effet, quelles que soient les parties mesurables $A$ et $B$ de $\{-1, 1\}^{-\mathbf{N}}$,

$$P[\xi_{n]} \in A; \eta_{n]} \in B; \zeta_{n+1} = 1 | \sigma(\zeta)] = \mathbf{I}_{[\zeta_{n+1}=1]} \nu_n(\zeta, A) \nu_n(\zeta, B)$$

par indépendance conditionnelle de $\xi$ et $\eta$ sachant $\zeta$. Mais l'étape 1 de la démonstration du Lemme 1 montre que $\nu_n(\zeta, A) \nu_n(\zeta, B)$ est une fonction de $\zeta_{n]}$, donc est indépendante de $\zeta_{n+1}$. En passant aux espérances, on obtient donc

$$P[\xi_{n]} \in A; \eta_{n]} \in B; \zeta_{n+1} = 1] = \tfrac{1}{2} \mathbf{E}[\nu_n(\zeta, A) \nu_n(\zeta, B)]$$
$$= \tfrac{1}{2} P[\xi_{n]} \in A; \eta_{n]} \in B].$$

D'après la Proposition 5, la suite $(X_n, Y_n)_{n \in \mathbf{Z}}$ définie par $X_n = \xi_{n]}, Y_n = \eta_{n]}$ est une chaîne de Markov sur $E^2$ pour les noyaux de transition $(N_n)_{n \in \mathbf{Z}}$. L'hypothèse (3) entraîne que les suites $(X_n)_{n \in \mathbf{Z}}$ et $(Y_n)_{n \in \mathbf{Z}}$ sont égales presque sûrement, donc que $\xi = \eta$ presque sûrement. Par conséquent, $\nu(\alpha, \cdot)$ est une masse de Dirac pour $\mu$-presque tout $\alpha$. Donc il existe une application mesurable $F$ de $\{-1; 1\}^{\mathbf{Z}}$ dans $\{-1; 1\}^{\mathbf{Z}}$ telle que $\varepsilon = F(H\varepsilon) = F(\Phi(\varepsilon))$ presque sûrement, ce qui entraîne l'égalité $\mathcal{F}_{\infty}^{H\varepsilon} = \mathcal{F}_{\infty}^{\varepsilon}$.  $\square$

### 1.5. *Cas d'une transformation homogène.*

Gardons les notations du paragraphe précédent. On dit que la transformation $\Phi$ est homogène lorsque les applications $h_n$ sont toutes égales à une même application $h$. Dans ce cas, les applications $f_+^{(n)}$, $f_-^{(n)}$ et le noyau de transition $N_n$ ne dépendent pas de $n$ et on les note simplement $f_+$, $f_-$ et $N$.

PROPOSITION 7. *Si la transformation $\Phi$ est homogène, on peut remplacer la condition* (3) *du Théorème* 6 *par la condition suivante: "toute chaîne de Markov stationnaire indexée par $\mathbf{Z}$ pour le noyau de transition $N$ vit sur la diagonale de $E^2$ (presque sûrement)."*

DÉMONSTRATION. Il suffit de vérifier que la démonstration de l'implication (3) $\Rightarrow$ (1) reste valable lorsqu'on remplace la condition (3) par cette condition apparemment plus faible. Reprenons donc trois jeux de pile ou face $\xi, \eta, \zeta$ tels que conditionnellement à $\zeta = \alpha$, les variables aléatoires $\xi$ et $\eta$ soient indépendantes de loi $\nu(\alpha, \cdot)$. Alors la loi des couples $(\xi, \zeta)$ et $(\eta, \zeta)$ est égale à la loi du couple $(\varepsilon, \Phi(\varepsilon))$, qui est invariante par translation temporelle, puisque la transformation $\Phi$ est homogène. Si $T$ est une translation temporelle sur $\{-1, 1\}^{\mathbf{Z}}$, alors $\sigma(\zeta) = \sigma(T(\zeta))$, donc l'indépendance conditionnelle de $\xi$ et $\eta$ sachant $\zeta$ entraîne l'indépendance conditionnelle de $T(\xi)$ et $T(\eta)$ sachant $T(\zeta)$. On en déduit facilement que la loi du triplet $(\xi, \eta, \zeta)$ est elle aussi invariante par translation temporelle. Par conséquent, la suite $(X_n, Y_n)_{n \in \mathbf{Z}}$ définie par $X_n = \xi_{n]}$, $Y_n = \eta_{n]}$ est une chaîne de Markov stationnaire sur $E^2$ pour le noyau de transition $N$. On peut donc affaiblir l'hypothèse (3) par la restriction "stationnaire." $\square$

L'intérêt de cette condition affaiblie est qu'elle se reformule simplement en terme de probabilités invariantes de la chaîne. Par ailleurs, l'énoncé devient plus simple si on transporte la chaîne de Markov de $E^2$ à $[0,1]^2$ ou à $[0,1[^2$ à l'aide des développements dyadiques.

En effet, notons $E'$ l'ensemble des suites $x \in \{-1, 1\}^{-\mathbf{N}}$ telles que $x_n = -1$ pour une infinité d'entiers $n$. Comme tout jeu de pile ou face indexé par $-\mathbf{N}$ appartient presque sûrement à $E'$, on peut remplacer $E$ par $E'$ dans l'énoncé de la Proposition 7.

Soit $S$ la bijection de $E'$ dans $[0,1[$ définie par

$$S(x) = \sum_{n=-\infty}^{0} 2^{n-1} \frac{1 + x_n}{2} \qquad \text{pour tout } x = (x_n)_{n \in -\mathbf{N}} \in \{-1, 1\}^{-\mathbf{N}}$$

et $\xi : [0,1[ \to E'$ la bijection réciproque de $S$. Notons enfin $A^+ = \{t \in [0,1[ : h(\xi(t)) = 1\}$ et $A^- = \{t \in [0,1[ : h(\xi(t)) = -1\}$. Alors pour tout $x \in E'$ et pour tout $\alpha \in \{+, -\}$,

$$S(f_\alpha(x)) = \sum_{n=-\infty}^{-1} 2^{n-1} \frac{1 + x_{n+1}}{2} + 2^{-1} \frac{1 + \alpha h(x)}{2} = \frac{S(x)}{2} + \frac{\mathbf{I}_{A^\alpha}(S(x))}{2}.$$

Cette égalité montre que si $(X_n, Y_n)_{n \in \mathbf{Z}}$ est une chaîne de Markov stationnaire sur $E'^2$ pour le noyau de transition $N$, alors $(S(X_n), S(Y_n))_{n \in \mathbf{Z}}$ est une chaîne de Markov stationnaire sur $[0,1[^2$ pour le noyau de transition $M$ défini par

$$M((u,v), \cdot) = \tfrac{1}{2}(\delta_{((u+\mathbf{I}_{A^+}(u)/2, (v+\mathbf{I}_{A^+}(v))/2)} + \delta_{((u+\mathbf{I}_{A^-}(u))/2, (v+\mathbf{I}_{A^-}(v))/2)}).$$

Remarquons qu'un pas pour cette chaîne correspond à une homothétie de rapport $\frac{1}{2}$, dont le centre est choisi au hasard parmi $(0,0)$ et $(1,1)$ si $(u,v) \in (A^+ \times A^+) \cup (A^- \times A^-)$, parmi $(0,1)$ et $(1,0)$ si $(u,v) \in (A^+ \times A^-) \cup (A^- \times A^+)$.

Réciproquement, si $(U_n, V_n)_{n \in \mathbf{Z}}$ est une chaîne de Markov stationnaire sur $[0,1[^2$ pour le noyau de transition $M$, alors $(\xi(U_n), \xi(V_n))_{n \in \mathbf{Z}}$ est une chaîne de Markov stationnaire sur $E'^2$ pour le noyau de transition $N$. On en déduit la caractérisation suivante:

PROPOSITION 8. *Pour que $\mathcal{F}^{H\varepsilon}_\infty = \mathcal{F}^\varepsilon_\infty$, il faut et il suffit que la seule probabilité invariante pour le noyau de transition $M$ soit la loi uniforme sur la diagonale de $[0,1[^2$.*

*Retour sur l'exemple du Paragraphe* 1.2. Reprenons les notations du Paragraphe 1.2. Alors tout $t \in [0,1[$, la suite $\xi(t)$ n'est autre que la suite $(\varepsilon_n(t))_{n \leq 0}$. Par conséquent,

$$\sum_{n=-\infty}^{-1} 2^{n-1} \frac{1 + \varepsilon_n(t)}{2} = S(\xi(t)) = t.$$

D'après l'expression de $h$ obtenue à la fin du paragraphe 1.2, on a donc

$$h(\xi(t)) = \begin{vmatrix} 1, & \text{si } t \in [0, \frac{1}{3}[ \cap [\frac{2}{3}, 1[, \\ -1, & \text{si } t \in [\frac{1}{3}, \frac{2}{3}[. \end{vmatrix}$$

Autrement dit, $A^+ = [0, \frac{1}{3}[ \cap [\frac{2}{3}, 1[$ et $A^- = [\frac{1}{3}, \frac{2}{3}[$.

On voit donc (figure 1) que $(A^+ \times A^+) \cup (A^- \times A^-)$ est une réunion de 5 carrés de côté $\frac{1}{3}$ obtenus en divisant le carré $[0,1]^2$ en 9, et $(A^+ \times A^-) \cup (A^- \times A^+)$ la réunion des 4 autres. Figure 1 qui suit montre le comportement d'une chaîne de Markov de noyau $M$ sur $[0,1]^2$. Si la chaîne a pour loi initiale la loi uniforme sur $[0,1]^2$, la loi après un pas est la loi uniforme sur une réunion de 18 carrés de côté $\frac{1}{6}$, la loi après deux pas est la loi uniforme sur une réunion de 36 carrés de côté $\frac{1}{12}$.... On voit que la loi après $n$ pas tend vers la loi uniforme sur $\{(u,v) \in [0,1]^2 : v - u \in \frac{1}{3}\mathbf{Z}\}$, qui est une probabilité invariante pour le noyau de transition $M$.

## 2. Transformations homogènes de longueur finie: principaux résultats.

Dans cette partie, on fixe $d \in \mathbf{N}^*$ et une application $\phi$ de $\{-1;1\}^d$ dans $\{-1;1\}$. On note $\Phi : \{-1;1\}^{\mathbf{Z}} \to \{-1;1\}^{\mathbf{Z}}$ la transformation définie par

$$\Phi(x) = (\phi(x_{n-d}, \ldots, x_{n-1})x_n)_{n \in \mathbf{Z}} \qquad \text{pour tout } x = (x_n)_{n \in \mathbf{Z}} \in \{-1;1\}^{\mathbf{Z}}.$$

Pour tout $\alpha \in \{-1;1\}$, on note $f_\alpha$ l'application de $\{-1;1\}^d$ dans $\{-1;1\}^d$ définie par

$$f_\alpha(x_1, \ldots, x_d) = (x_2, \ldots, x_d, \alpha\phi(x_1, \ldots, x_d))$$

quel que soit $(x_1, \ldots, x_d) \in \{-1;1\}^d$. Par construction, les applications $f_+$ et $f_-$ vérifient

$$(\varepsilon_{n-d+1}, \ldots, \varepsilon_n) = f_{H_n \varepsilon_n}(\varepsilon_{n-d}, \ldots, \varepsilon_{n-1}),$$

avec $H_n = \phi(\varepsilon_{n-d}, \ldots, \varepsilon_{n-1})$.

2.1. *La relation d'accordabilité.* Adaptons l'énoncé du Théorème 6 à notre situation, dans laquelle $H_n$ ne dépend que de $(\varepsilon_{n-d}, \ldots, \varepsilon_{n-1})$ et les noyaux de transition $N_n$ ne dépendent pas de $n$.

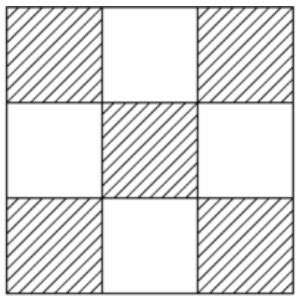

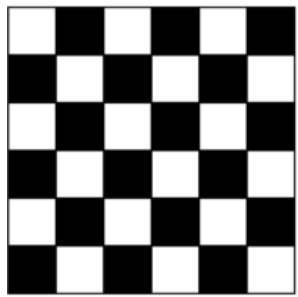

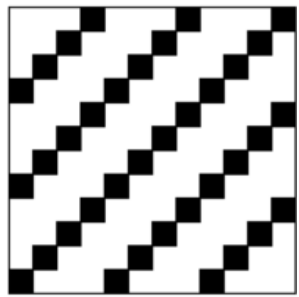

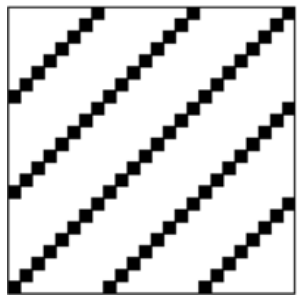

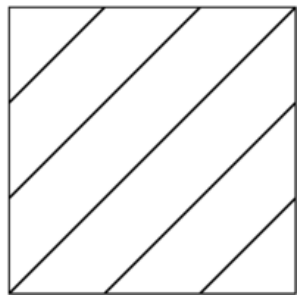

Fig. 1. *Le premier carré indique les parties* $(A^+ \times A^+) \cup (A^- \times A^-)$ *(carrés hachurés) et* $(A^+ \times A^-) \cup (A^- \times A^+)$ *(carrés blancs). Les suivants indiquent (en noir) les positions atteignables en un pas, en deux pas, en trois pas, et après une infinité de pas.*

PROPOSITION 9.   *Dans l'énoncé du Théorème 6, on peut remplacer l'espace $E = \{-1,1\}^{-\mathbf{N}}$ par $E_d = \{-1;1\}^d$ et les noyaux de transition $N_n$ par les probabilités de transition définies par*

$$p((a,b),(f_+(a),f_+(b))) = p((a,b),(f_-(a),f_-(b))) = \tfrac{1}{2}$$

*pour tout $(a,b) \in E_d^2$.*

En effet, la donnée d'une chaîne de Markov indexée par $\mathbf{Z}$ sur $E_d = \{-1;1\}^d$ pour ces probabilités de transition $p$ est équivalente à la donnée d'une chaîne de Markov indexée par $\mathbf{Z}$ sur $E = \{-1,1\}^{-\mathbf{N}}$ pour les noyaux $N_n$.

L'avantage de cette formulation est de se ramener à l'étude d'une chaîne de Markov sur un espace fini, ce qui nous amène à définir la relation d'accordabilité.

DÉFINITION 10.   Si $a,b \in E_d$, on dit que $a$ et $b$ sont accordables si et seulement si il existe une suite finie de signes $\alpha_1, \ldots, \alpha_l$ tels que la composée $f_{\alpha_l} \circ \cdots \circ f_{\alpha_1}$ prenne la même valeur en $a$ et $b$. Plus généralement on dit que des éléments $a_1, \ldots, a_m \in E_d$ sont simultanément accordables si et seulement si il existe une suite finie de signes $\alpha_1, \ldots, \alpha_l$ tels que la composée $f_{\alpha_l} \circ \cdots \circ f_{\alpha_1}$ prenne la même valeur en $a_1, \ldots, a_m$.

L'accordabilité simultanée de tous les éléments de $E_d$ coïncide avec la notion de "point collapsing semigroup" (pour le semigroupe associé à $f_+$ et $f_-$) introduite par Rosenblatt [5].

*Attention*!   La relation d'accordabilité n'est pas transitive en général. L'accordabilité simultanée est donc plus forte que l'accordabilité deux-à-deux.

En termes de chaînes de Markov sur $E_d^2$, dire que $a$ et $b$ sont accordables revient à dire que la diagonale de $E_d^2$ est atteignable depuis $(a,b)$. Comme la diagonale est une classe récurrente, tout couple $(a,b)$ formé de deux éléments accordables distincts est un état transitoire.

Ces remarques permettent d'obtenir une condition nécessaire et suffisante de conservation de l'information plus opératoire que celle de la Proposition 9. Le théorème ci-dessous, du moins l'équivalence (1) ⇔ (3), peut être vu comme une conséquence du Théorème 1 de Rosenblatt [5], même si la démonstration de l'implication (1) ⇒ (3) y est pour le moins elliptique. Nous faisons le point sur les liens entre le Théorème 1 de Rosenblatt [5] et nos résultats dans le paragraphe suivant.

THÉORÈME 11.   *Il y a équivalence entre:*

(1) $\mathcal{F}_\infty^{H\varepsilon} = \mathcal{F}_\infty^\varepsilon$;

(2) *tous les éléments de $E_d$ sont deux-à-deux accordables;*

(3) *tous les éléments de $E_d$ sont simultanément accordables.*

DÉMONSTRATION.   L'implication $(1) \Rightarrow (2)$ est une conséquence de la Proposition 9: si toute chaîne de Markov indexée par $\mathbf{Z}$ pour les probabilités de transition $p(\cdot, \cdot)$ vit sur la diagonale, la diagonale est la seule classe récurrente. Par conséquent, toute chaîne de Markov (indexée par $\mathbf{N}$) partant de $(a, b) \in E_d^2$ atteint en un temps fini la diagonale, ce qui montre que $a$ et $b$ sont accordables.

On prouve l'implication $(2) \Rightarrow (3)$ en démontrant par récurrence sur $m \in [2 \ldots 2^d]$ que $m$ points quelconques $a_1, \ldots, a_m \in E_d$ sont simultanément accordables. En effet, si $g = f_{\alpha_l} \circ \cdots \circ f_{\alpha_1}$ accorde $a_{m-1}$ et $a_m$, alors il y a au plus $m - 1$ points distincts parmi $g(a_1), \ldots, g(a_m)$. L'accordabilité simultanée de $g(a_1), \ldots, g(a_m)$ entraîne celle de $a_1, \ldots, a_m$.

L'implication $(3) \Rightarrow (1)$ est une conséquence immédiate de la Proposition 9 et du fait qu'une chaîne de Markov indexée par $\mathbf{Z}$ sur un espace d'états fini vit sur la réunion des classes récurrentes.   $\square$

*Une démonstration directe et constructive de l'implication* $(3) \Rightarrow (1)$. Voyons une autre démonstration indiquant comment retrouver la suite $(\varepsilon_n)_{n \in \mathbf{Z}}$ à partir de la suite $(H_n \varepsilon_n)_{n \in \mathbf{Z}}$. Cette démonstration est semblable à celle de Rosenblatt [5].

Par hypothèse, on peut trouver une suite finie de signes $\alpha_1, \ldots, \alpha_l$ tels que la composée $f_{\alpha_l} \circ \cdots \circ f_{\alpha_1}$ soit constante. Comme $(H_n \varepsilon_n)_{n \in \mathbf{Z}}$ est un jeu de pile ou face, il existe presque sûrement une infinité d'instants $n \leq 0$ tels que $(H_{n-l+1} \varepsilon_{n-l+1}, \ldots, H_n \varepsilon_n) = (\alpha_1, \ldots, \alpha_l)$. Cela impose la valeur du $d$-uplet $(\varepsilon_{n-d+1}, \ldots, \varepsilon_n)$ à ces instants, et donc à tout instant, en vertu de la relation de récurrence $(\varepsilon_{n-d+1}, \ldots, \varepsilon_n) = f_{H_n \varepsilon_n}(\varepsilon_{n-d}, \ldots, \varepsilon_{n-1})$.

Plus précisément, si $f_{\alpha_l} \circ \cdots \circ f_{\alpha_1}$ est constante égale à $a$, alors pour tout $n \in \mathbf{Z}$,

$$(\varepsilon_{n-d+1}, \ldots, \varepsilon_n) = f_{H_n \varepsilon_n} \circ \cdots \circ f_{H_{k+1} \varepsilon_{k+1}}(a),$$

d'où

$$\varepsilon_n = H_n \varepsilon_n \phi(f_{H_{n-1} \varepsilon_{n-1}} \circ \cdots \circ f_{H_{k+1} \varepsilon_{k+1}}(a)),$$

où $k$ est le dernier instant avant $n$ tel que $(H_{k-l+1} \varepsilon_{k-l+1}, \ldots, H_k \varepsilon_k) = (\alpha_1, \ldots, \alpha_l)$. Cette égalité exprime la suite $(\varepsilon_n)_{n \in \mathbf{Z}}$ comme l'image de la suite $\Phi((\varepsilon_n)_{n \in \mathbf{Z}}) = (H_n \varepsilon_n)_{n \in \mathbf{Z}}$ par une transformation homogène $\Psi$. En notant $\mu$ la loi uniforme sur $\{-1, 1\}^{\mathbf{Z}}$, on a donc $\Psi \circ \Phi = \mathrm{Id}$ $\mu$-presque partout. Mais comme la transformation $\Phi$ préserve la loi $\mu$, on a aussi $\Phi \circ \Psi = \mathrm{Id}$ $\mu$-presque partout. Les transformations $\Phi$ et $\Psi$ sont donc inverses l'une de l'autre.

Exemple: prenons $d = 2$ et $\phi : (x_1, x_2) \mapsto \max(x_1, x_2)$, d'où $H_n = \max(\varepsilon_{n-2}, \varepsilon_{n-1})$ pour tout $n \in \mathbf{Z}$. La figure 2 montre les graphes de $f_+$ et $f_-$.

On voit que $f_+^2 \circ f_- \circ f_+$ envoie tous les points de $\{-1, 1\}^2$ sur $(1, 1)$. Par conséquent, la filtration naturelle de $(H_n \varepsilon_n)_{n \in \mathbf{Z}}$ est égale à celle de $(\varepsilon_n)_{n \in \mathbf{Z}}$. Pour reconstituer la suite $(\varepsilon_n)_{n \in \mathbf{Z}}$ connaissant $(H_n \varepsilon_n)_{n \in \mathbf{Z}}$, on cherche dans la suite $(H_n \varepsilon_n)_{n \in \mathbf{Z}}$ les instants où les quatre derniers bits sont $1, -1, 1, 1$. A ces instants, les deux derniers bits de la suite $(\varepsilon_n)_{n \in \mathbf{Z}}$ sont $1, 1$, ce qui permet de reconstituer toute la suite après ces instants.

2.2. *Précisions sur les liens avec les travaux de M. Rosenblatt.* Dans le Paragraphe 3 de [5], Rosenblatt s'intéresse aux chaînes de Markov uniformes sur un espace d'états dénombrable, c'est-à-dire aux chaînes de Markov dont les probabilités de transition issues d'un état et rangées par ordre décroissant ne dépendent pas de l'état considéré (autrement dit, les lignes de la matrice de transition se déduisent les unes des autres par permutation de leur coefficients).

Notons $E = \{1, 2, \ldots\}$ l'espace d'états et $(p_{i,j})_{i,j \in E}$ la matrice de transition. L'hypothèse d'uniformité se traduit par le fait que pour chaque état $i$, la liste ordonnée des valeurs $p_{i,j}$ est une suite $q_1 \geq q_2 \geq \cdots$ qui ne dépend pas de $i$.

Rosenblatt considère une chaîne de Markov stationnaire $(X_n)_{n \in \mathbf{Z}}$ de matrice de transition $(p_{i,j})_{i,j \in E}$. Il suppose en outre que les valeurs $q_k$ non nulles sont toutes différentes. Cela lui permet de construire (presque sûrement) une suite $(\xi_n)_{n \in \mathbf{Z}}$ de variables aléatoires indépendantes de même loi $\sum q_k \delta_k$ en posant pour tout $n \in \mathbf{Z}$ et pour tout couple $(i, j)$ tel que $p_{i,j} = q_k > 0$

$$\xi_n = k \qquad \text{sur l'événement } [X_{n-1} = i; X_n = j].$$

Les variables aléatoires $(\xi_n)_{n \in \mathbf{Z}}$ ainsi obtenues sont des "innovations" de la filtration $(\mathcal{F}_n^X)_{n \in \mathbf{Z}}$ au sens où pour tout $n \in \mathbf{Z}$, $\xi_n$ est indépendante de $\mathcal{F}_{n-1}^X$ et $\mathcal{F}_n^X = \mathcal{F}_{n-1}^X \vee \sigma(\xi_n)$. De plus, la suite $(X_n)_{n \in \mathbf{Z}}$ vérifie une relation de récurrence de la forme $X_n = f(X_{n-1}, \xi_n) = f_{\xi_n}(X_{n-1})$.

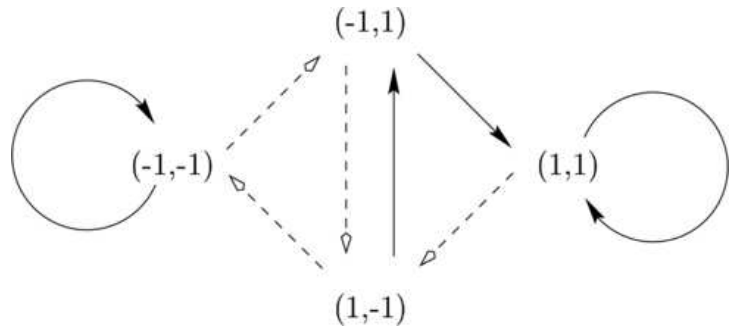

Fig. 2. *Les flèches pleines indiquent l'image par $f_+$, les flèches creuses l'image par $f_-$.*

Une question naturelle est alors de savoir si la suite d'innovations $(\xi_n)_{n \in \mathbf{Z}}$ engendre la même filtration que la chaîne de Markov $(X_n)_{n \in \mathbf{Z}}$. Dans son Théorème 1, Rosenblatt affirme qu'une condition nécessaire et suffisante est que le semigroupe engendré par les applications $f_k : x \mapsto f(x, k)$ soit "point collapsing," c'est-à-dire qu'il existe une composée d'applications $f_k$ qui est constante. En réalité, Rosenblatt regarde plutôt le semigroupe engendré par les matrices $M_k = (\mathbf{I}_{[p_{i,j} = q_k]})_{i,j \in E} = (\mathbf{I}_{[j = f(i,k)]})_{i,j \in E}$, mais cela revient au même.

Notre Théorème 11 peut être vu comme un cas particulier du théorème 1 de Rosenblatt, avec $E = \{-1, 1\}^d$, $q_1 = q_2 = \frac{1}{2}$, $f_1 = f_+$ et $f_2 = f_-$, $X_n = (\varepsilon_{n-d+1}, \ldots, \varepsilon_n)$, $\xi_n = 1$ ou 2 selon que $H_n \varepsilon_n = 1$ ou $-1$. Bien sûr, l'hypothèse que les valeurs $q_k$ non nulles sont toutes différentes n'est pas vérifiée dans notre situation, mais elle ne joue aucun rôle dans le Théorème 1 de Rosenblatt: elle a seulement servi à définir les variables $\xi_n$ sans avoir à choisir un ordre entre deux transitions de même probabilité partant d'un même état.

Cependant, Rosenblatt ne démontre vraiment qu'une implication: s'il existe une composée d'applications $f_k$ qui est constante, alors la suite $(\xi_n)_{n \in \mathbf{Z}}$ engendre la même filtration que la chaîne de Markov $(X_n)_{n \in \mathbf{Z}}$. Mais pour la réciproque, l'argument de Rosenblatt revient à admettre l'idée intuitive suivante: si l'on connaît la suite $(\xi_n)_{n \in \mathbf{Z}}$, la seule chose qu'on puisse dire de $X_n$ est que $X_n$ appartient à l'intersection des images des composées $f_{\xi_n} \circ f_{\xi_{n-1}} \circ \cdots \circ f_{\xi_{n-l+1}}$ pour $l \in \mathbf{N}$. Dans notre cas, le Théorème 20 donne un sens précis à cette affirmation et la démontre.

Pour démontrer l'implication $(1) \Rightarrow (3)$ dans le théorème 11 (qui est un cas particulier de l'implication non démontrée de Rosenblatt), nous avons utilisé un argument de couplage, en considérant une chaîne de Markov sur $E_d^2$ dont les deux composantes sont des chaînes de Markov sur $E_d$ gouvernées par les mêmes innovations. Dans la Proposition 9.3.3 de sa thèse [2], S. Laurent démontre l'implication non démontrée par Rosenblatt pour une chaîne de Markov uniforme et stationnaire sur un espace d'états fini. L'argument de Laurent (obtenu indépendamment de nous) est un argument de couplage différent du nôtre, utilisant un critère qu'il appelle critère de Doeblin.

Signalons enfin que l'obtention d'applications constantes par composition d'applications de type $f_{\xi_n}$ est une idée utilisée également dans d'autres buts: c'est par exemple le coeur de l'algorithme de couplage depuis le passé de Propp et Wilson [4] qui fournit une simulation "exacte" de la probabilité invariante d'une chaîne de Markov sur un espace d'états fini.

### 2.3. *L'algorithme de Roland Bacher.*
L'algorithme de Roland Bacher calcule la matrice d'adjacence de la relation d'accordabilité. Le principe con-siste à calculer par récurrence la matrice $(A_n[i,j])$ où $A_n[i,j]$ vaut 1 si $i$ et $j$ sont accordables en au plus $n$ étapes et 0 sinon. On part de $A_0 = I$ et on passe de $A_n$ à $A_{n+1}$ en changeant $A_n[i,j]$ de 0 en 1 si et seulement si $A_n[i,j] = 0$ et si $A_n[f_+(i), f_+(j)] = 1$ ou $A_n[f_-(i), f_-(j)] = 1$. On s'arrête dès que $A_n = A_{n+1}$.

REMARQUE. En fait, l'algorithme de Roland Bacher est une variante plus rapide et plus simple à programmer de l'algorithme ci-dessus, dans laquelle les coefficients sont modifiés de façon dynamique.

LEMME 12. *Deux points accordables sont accordables en au plus* $2^{d-1} \times (2^d - 1)$ *étapes.*

DÉMONSTRATION. Il suffit de remarquer que le nombre de coefficients égaux à 1 dans la matrice $A_n$ varie entre $2^d$ et $2^{2d}$ et qu'il croît strictement par multiples de 2 avant d'atteindre la valeur limite. Un autre argument est que si $a$ et $b$ sont accordables, alors un chemin minimal de $(a,b)$ à la diagonale pour la chaîne de Markov sur $E_d^2$ visite au plus une fois chaque paire de la forme $\{(c,d); (d,c)\}$. □

La borne ci-dessus est loin d'être optimale, et on peut l'améliorer facilement en utilisant la structure du graphe associé aux applications $f_+$ et $f_-$. Par définition de $f_+$ et de $f_-$, tout élément $(y_1, \ldots, y_d) \in E_d$ possède exactement deux antécédents par $f_+$ et $f_-$ réunis: ce sont les $d$-uplets $(1, y_1, \ldots, y_{d-1})$ et $(-1, y_1, \ldots, y_{d-1})$. De plus, les antécédents de $(y_1, \ldots, y_{d-1}, y_d)$ par $f_+$ et $f_-$ sont aussi les antécédents de $(y_1, \ldots, y_{d-1}, -y_d)$ par $f_-$ et $f_+$. Cela nous amène à poser la définition suivante.

DÉFINITION 13 (Points de convergence).

(1) On appelle point de convergence positif tout élément de $E_d$ ayant deux antécédents par $f_+$ (et donc aucun pour $f_-$).

(2) On appelle point de convergence négatif tout élément de $E_d$ ayant deux antécédents par $f_-$ (et donc aucun pour $f_+$).

(3) On appelle point ordinaire tout élément de $E_d$ ayant un antécédent par $f_+$ et un par $f_-$.

PROPOSITION 14. *Deux points accordables sont accordables en au plus* $2^{2d-2}$ *étapes.*

DÉMONSTRATION. Notons $S$ l'ensemble des paires d'éléments de $E_d$. On vérifie immédiatement que pour tout $\alpha \in \{-1, 1\}$ et $x, y \in E_d$, l'égalité $f_+(x) = f_\alpha(y)$ équivaut à $f_-(x) = f_{-\alpha}(y)$. Cette remarque permet de définir une relation d'équivalence sur $S$ par

$$\{a, b\} \sim \{c, d\} \iff \exists \alpha \in \{-1, 1\} : \{f_+(a), f_+(b)\} = \{f_\alpha(c), f_\alpha(d)\}.$$

Il est clair que si $a$ et $b$ sont accordables, alors un chemin minimal de $(a, b)$ à la diagonale pour la chaîne de Markov sur $E_d^2$ visite au plus une fois une classe d'équivalence. Le nombre minimum d'étapes est donc majoré par le nombre de classes d'équivalence.

Etudions la classe d'équivalence d'une paire $\{a, b\} \in S$ quelconque, suivant que $f_+(a)$ et $f_+(b)$ sont ou non des points de convergence. Pour tout $x = (x_1, \ldots, x_d) \in E_d$, on notera $\widetilde{x} = (-x_1, x_2, \ldots, x_d)$. Soit $N$ le nombre de points de convergence positifs. Alors $N$ est aussi le nombre de points de convergence négatifs, donc il y a exactement $2^d - 2N$ points ordinaires.

Si $f_+(a) = f_+(b)$, alors ce point est un point de convergence et la classe de $\{a, b\}$ est réduite à $\{a, b\}$. Il y a $N$ classes de cette forme (soit $N$ paires).

Si $f_+(a) \neq f_+(b)$ et si $f_+(a)$ et $f_+(b)$ sont des points de convergence, alors la classe de $\{a, b\}$ contient quatre éléments: $\{a, b\}$, $\{\widetilde{a}, b\}$ $\{a, \widetilde{b}\}$ et $\{\widetilde{a}, \widetilde{b}\}$. Il y a $N(N-1)/2$ classes de cette forme [soit $2N(N-1)$ paires].

Si $f_+(a) \neq f_+(b)$ et s'il n'y a qu'un point de convergence parmi $f_+(a)$ et $f_+(b)$, alors la classe de $\{a, b\}$ contient deux éléments: $\{a, b\}$ et $\{\widetilde{a}, b\}$ si $f_+(a)$ est le point de convergence et $f_+(b)$ le point ordinaire. Il y a $N(2^d - 2N)$ classes de cette forme [soit $2N(2^d - 2N)$ paires].

Si $f_+(a) \neq f_+(b)$ et si $f_+(a)$ sont des points ordinaires, alors la classe de $\{a, b\}$ est réduite à $\{a, b\}$ et $\{\widetilde{a}, \widetilde{b}\}$. Mais ces paires sont confondues lorsque $f_+(a) = f_-(b)$ et distinctes dans le cas contraire. On obtient donc $(2^d - 2N)/2$ classes à un élément [soit $(2^d - 2N)/2$ paires] et $(2^d - 2N)(2^d - 2N - 2)/4$ classes à deux éléments [soit $(2^d - 2N)(2^d - 2N)/2$ paires].

On vérifie bien que le nombre total de paires est

$$N + 2N(N-1) + 2N(2^d - 2N) + \frac{2^d - 2N}{2} + \frac{(2^d - 2N)(2^d - 2N - 2)}{2},$$

soit $\frac{2^d(2^d - 1)}{2}$. Le nombre total de classes est

$$N + \frac{N(N-1)}{2} + N(2^d - 2N) + \frac{2^d - 2N}{2} + \frac{(2^d - 2N)(2^d - 2N - 2)}{4},$$

soit $\frac{4^d - 2N^2 + 2N}{4}$. Il est donc majoré par $4^{d-1}$.  $\square$

REMARQUE. La borne $2^{2d-2}$ sur le nombre d'étapes nécessaires pour accorder deux éléments de $E_d$ est optimale si $d = 2$. En effet, dans l'exemple du Paragraphe 2.1, on voit que $(1, -1)$ et $(-1, -1)$ ne peuvent pas s'accorder en moins de quatre coups. Pour $d = 3$ et $\phi : (x_1, x_2, x_3) \mapsto \max(x_1, -x_2)$, neuf coups sont nécessaires pour accorder les triplets $(1, 1, -1)$ et $(-1, 1, -1)$.

2.4. *Etude d'une chaîne de Markov sur l'ensemble des parties de $E_d$.* Pour décrire l'information perdue par la transformation $\Phi$, nous allons utiliser une chaîne de Markov sur l'ensemble des parties de $E_d$, noté $\mathcal{P}(E_d)$.

Soit $(\eta_n)_{n \in \mathbf{N}^*}$ une suite de variables indépendantes de loi $(\delta_{-1} + \delta_1)/2$. Pour $A \in \mathcal{P}(E_d)$, notons $(X_n^A)_{n \in \mathbf{N}}$ la suite de variables aléatoires définie par $X_0^A = A$ et la relation de récurrence $X_n^A = f_{\eta_n}(X_{n-1}^A)$. Alors $(X_n^A)_{n \in \mathbf{N}}$ est une chaîne de Markov sur $\mathcal{P}(E_d)$ pour les probabilités de transition définies par

$$q(A, f_+(A)) = q(A, f_-(A)) = \tfrac{1}{2}.$$

Bien entendu, la suite des cardinaux des parties $(X_n^A)_{n \in \mathbf{N}}$ est décroissante. Par ailleurs, si deux parties de $E_d$ vérifient $A \subset B$, alors $X_n^A \subset X_n^B$ pour tout $n \in \mathbf{N}$.

LEMME 15. *Soit $M$ le nombre maximum d'éléments deux à deux non accordables. Alors presque sûrement, $\mathrm{Card}(X_n^{E_d}) = M$ à partir d'un certain rang.*

DÉMONSTRATION. La démonstration comporte deux parties:

Si une partie $A$ ne contient que des points deux-à-deux non accordables, alors pour toute suite finie de signes $\alpha_1, \ldots, \alpha_n$, les images par $f_{\alpha_n} \circ \cdots \circ f_{\alpha_1}$ des éléments de $A$ sont deux-à-deux distinctes, donc $\mathrm{Card}(X_n^A) = \mathrm{Card}(A)$ pour tout $n \in \mathbf{N}$. En prenant pour $A$ un ensemble de $M$ points deux-à-deux non accordables, on obtient $\mathrm{Card}(X_n^{E_d}) \geq \mathrm{Card}(X_n^A) = M$.

Si une partie $A$ est un état récurrent de la chaîne, alors la décroissance de $\mathrm{Card}(X_n^A)$ impose que $\mathrm{Card}(X_n^A) = \mathrm{Card}(A)$ pour tout $n \in \mathbf{N}$. Par conséquent, $A$ ne contient que des points deux-à-deux non accordables (sinon, le cardinal diminuerait avec probabilité strictement positive), d'où $\mathrm{Card}(A) \leq M$. Comme $(X_n^{E_d})_{n \in \mathbf{N}}$ est une chaîne de Markov sur un espace d'états fini, elle finit presque sûrement dans une composante récurrente d'où $\mathrm{Card}(X_n^{E_d}) \leq M$ à partir d'un certain rang. $\square$

En choisissant un entier naturel $l$ tel que $P[\mathrm{Card}(X_l^{E_d}) = M] > 0$, on obtient:

COROLLAIRE 16. *Soit $M$ le nombre maximum d'éléments deux à deux non accordables. Il existe une suite finie de signes $\alpha_1, \ldots, \alpha_l$ telle que $f_{\alpha_l} \circ \cdots \circ f_{\alpha_1}(E_d)$ soit un ensemble formé de $M$ éléments deux-à-deux non accordables.*

COROLLAIRE 17. *Soit $M$ le nombre maximum d'éléments deux à deux non accordables et $\mathcal{C}$ l'ensemble des états récurrents de $\mathcal{P}(E_d)$ à $M$ éléments. Alors $\mathcal{C}$ est une classe récurrente.*

DÉMONSTRATION.    En effet, soient $\alpha_1, \ldots, \alpha_l \in \{-1, 1\}$ tels que l'ensemble $A = f_{\alpha_l} \circ \cdots \circ f_{\alpha_1}(E_d)$ soit formé de $M$ éléments deux-à-deux non accordables. Alors pour tout $B \in \mathcal{C}$, $f_{\alpha_l} \circ \cdots \circ f_{\alpha_1}(B) \subset A$ et Card$(f_{\alpha_l} \circ \cdots \circ f_{\alpha_1}(B)) = M$, d'où $f_{\alpha_l} \circ \cdots \circ f_{\alpha_1}(B) = A$, ce qui montre que l'état $A$ est accessible depuis tout état $B \in \mathcal{C}$. Par conséquent, $A$ communique avec tout état $B \in \mathcal{C}$.    □

Nous allons maintenant nous intéresser aux chaînes de Markov indexées par $\mathbf{Z}$ sur $\mathcal{C}$ (pour les probabilités de transition définies au début du paragraphe):

PROPOSITION 18.    *Soit $(X_n)_{n \in \mathbf{Z}}$ une chaîne de Markov sur $\mathcal{C}$, gouvernée par un jeu de pile ou face $(\eta_n)_{n \in \mathbf{Z}}$: pour tout $n \in \mathbf{Z}$, $X_n = f_{\eta_n}(X_{n-1})$. Alors il existe une application mesurable $F: \{-1, 1\}^{-\mathbf{N}} \to \mathcal{C}$ telle que presque sûrement $X_n = F(\eta_{n]})$ pour tout $n \in \mathbf{Z}$.*

DÉMONSTRATION.    Soit $\alpha_1, \ldots, \alpha_l$ une suite finie de signes comme dans le Corollaire 16 et $A = f_{\alpha_l} \circ \cdots \circ f_{\alpha_1}(E_d)$. Alors presque sûrement $(\eta_{k-l+1}, \ldots, \eta_k) = (\alpha_1, \ldots, \alpha_l)$ pour une infinité de $k \in -\mathbf{N}$, d'où $X_k = A$ à ces instants $k$. Pour tout $n \in \mathbf{Z}$ et pour tout $k \leq n$, on a donc $X_n = f_{\eta_n} \circ \cdots \circ f_{\eta_{k+1}}(A)$ où $k$ est le dernier instant avant $n$ tel que $(\eta_{k-l+1}, \ldots, \eta_k) = (\alpha_1, \ldots, \alpha_l)$. Cette écriture fournit une application $F$ (définie presque partout) telle que pour tout $n \in \mathbf{Z}$, $X_n = F(\eta_{n]})$ presque sûrement.    □

2.5.  *Description de l'information perdue.*    Nous allons voir que l'application $F$ fournie par la Proposition 18 joue un rôle-clé dans la description de l'information perdue. Énonçons déjà un premier résultat dans ce sens:

PROPOSITION 19.    *Soit $F: \{-1, 1\}^{-\mathbf{N}} \to \mathcal{C}$ l'application fournie par la Proposition 18. Alors presque sûrement, $(\varepsilon_{n-d+1}, \ldots, \varepsilon_n) \in F((H\varepsilon)_{n]})$ pour tout $n \in \mathbf{Z}$.*

DÉMONSTRATION.    La démonstration est semblable à celle de la Proposition 18, et utilise la relation de récurrence

$$(\varepsilon_{n-d+1}, \ldots, \varepsilon_n) = f_{H_n \varepsilon_n}(\varepsilon_{n-d}, \ldots, \varepsilon_{n-1}).$$

Presque sûrement $(H_{k-l+1}\varepsilon_{k-l+1}, \ldots, H_k \varepsilon_k) = (\alpha_1, \ldots, \alpha_l)$ pour une infinité de $k \in -\mathbf{N}$, d'où $(\varepsilon_{k-d+1}, \ldots, \varepsilon_k) \in A = F((H\varepsilon)_{k]})$ à ces instants $k$. On en déduit par récurrence le résultat pour les instants ultérieurs.    □

THÉORÈME 20.    *Soit $F: \{-1, 1\}^{-\mathbf{N}} \to \mathcal{C}$ l'application fournie par la Proposition 18. Alors presque sûrement, la loi de $(\varepsilon_{n-d+1}, \ldots, \varepsilon_n)$ sachant $\Phi(\varepsilon)$ est la loi uniforme sur $F((H\varepsilon)_{n]})$.*

Avant de démontrer le théorème, commençons par établir un lemme simple:

LEMME 21. *Pour tout $\alpha \in \mathbf{Z}$, $\Phi^{-1}(\alpha)$ est un ensemble fini de cardinal $\leq 2^d$.*

DÉMONSTRATION. En effet, soient $x^1, \ldots, x^{2^d+1}$ des éléments de $\Phi^{-1}(\alpha)$. Pour tout $n \geq 0$, il existe deux indices $i < j$ dans $[1 \ldots 2^d + 1]$ tels que $(x_n^i, \ldots, x_{n+d-1}^i) = (x_n^j, \ldots, x_{n+d-1}^j)$, puisque $\mathrm{Card}(E_d) = 2^d$. Choisissons deux tels indices et notons-les $i_n$ et $j_n$. Une récurrence immédiate montre que pour tout $n \geq 0$, les suites $x^{i_n}$ et $x^{j_n}$ coïncident après l'instant $n$. Mais la suite $(i_n, j_n)_{n \leq 0}$ est à valeurs dans un ensemble fini et prend donc une infinité de fois une même valeur $(i, j)$. Les suites $x^i$ et $x^j$ correspondantes sont donc égales. $\square$

Passons maintenant à la démonstration du Théorème 20:

DÉMONSTRATION DU THÉORÈME 20. Soit $M$ le nombre maximum d'éléments deux à deux non accordables. On sait d'après le Théorème 1 et le Lemme 21 que la loi de $\varepsilon$ sachant $\Phi(\varepsilon)$ est presque sûrement uniforme sur un ensemble fini de cardinal constant. Notons $N$ ce cardinal. La démonstration comprend trois étapes:

*Étape 1.* Pour tout $n \in \mathbf{Z}$, la loi de $(\varepsilon_k)_{k \geq n}$ sachant $\Phi(\varepsilon)$ est presque sûrement une loi uniforme sur un ensemble à $N$ éléments.

En effet, notons $N_n$ le nombre d'atomes de la loi de $(\varepsilon_k)_{k \geq n}$ sachant $\Phi(\varepsilon)$. Comme cette loi est (presque sûrement) l'image de la loi de $\varepsilon$ sachant $\Phi(\varepsilon)$ par la projection canonique de $\{-1, 1\}^{\mathbf{Z}}$ sur $\{-1, 1\}^{[n \ldots \infty[}$, on en déduit que $N_n \leq N$ presque sûrement, et il suffit de montrer que $N_n = N$ presque sûrement.

Comme la loi du couple $(\varepsilon, \Phi(\varepsilon))$ est invariante par translation temporelle, la loi de $N_n$ est indépendante de $n$. Et comme la suite $(N_n)_{n \in \mathbf{Z}}$ est presque sûrement décroissante [puisque la loi $(\varepsilon_k)_{k \geq n+1}$ sachant $\Phi(\varepsilon)$ est l'image de la loi de $(\varepsilon_k)_{k \geq n}$ sachant $\Phi(\varepsilon)$ par la projection canonique de $\{-1, 1\}^{[n \ldots \infty[}$ sur $\{-1, 1\}^{[n+1 \ldots \infty[}$], elle est presque sûrement constante.

Notons $a_1(\alpha), \ldots, a_N(\alpha)$ les atomes de la loi de $\varepsilon$ sachant que $\Phi(\varepsilon) = \alpha$, numérotés de sorte que les applications $a_k$ soient mesurables de $\{-1, 1\}^{\mathbf{Z}}$ dans $\{-1, 1\}^{\mathbf{Z}}$. Alors pour tout $n \in \mathbf{Z}$, l'événement $[N_n = N]$ est presque sûrement égal à l'événement "les restrictions de $a_1(\Phi(\varepsilon)), \ldots, a_N(\Phi(\varepsilon))$ à $[n \ldots \infty[$ sont deux-à-deux distinctes." Donc $P[N_n = N] \to 1$ quand $n \to -\infty$. Comme $N_n$ ne dépend pas de $n$ (presque sûrement), on a donc $N_n = N$ presque sûrement.

*Étape* 2. Pour tout $n \in \mathbf{Z}$, la loi de $(\varepsilon_n, \ldots, \varepsilon_{n+d-1})$ sachant $\Phi(\varepsilon)$ est presque sûrement une loi uniforme sur un ensemble à $N$ éléments.

Cela découle du fait que la loi de $(\varepsilon_n, \ldots, \varepsilon_{n+d-1})$ sachant $\Phi(\varepsilon)$ est presque sûrement l'image de la loi de $(\varepsilon_k)_{k \geq n}$ sachant $\Phi(\varepsilon)$ par l'application $(x_k)_{k \geq n} \mapsto (x_n, \ldots, x_{n+d-1})$ et que la restriction de cette application à $\Phi^{-1}(\alpha)$ est injective pour tout $\alpha \in \{-1, 1\}^{\mathbf{Z}}$.

*Étape* 3. Pour tout $n \in \mathbf{Z}$, la loi de $(\varepsilon_{n-d+1}, \ldots, \varepsilon_n)$ sachant $\Phi(\varepsilon)$ est la loi uniforme sur $F((H\varepsilon)_{n]})$.

D'après la Proposition 19, cette loi est portée (presque sûrement) par $F((H\varepsilon)_{n]})$, qui possède $M$ éléments. On a donc $N \leq M$ presque sûrement et il reste à montrer que $M \leq N$ presque sûrement, c'est-à-dire que la loi de $\varepsilon$ sachant $\Phi(\varepsilon)$ possède au moins $M$ atomes.

Pour cela, on considère une chaîne de Markov sur l'ensemble $S$ des $M$-uplets d'éléments de $E_d$ deux-à-deux non accordables. Si $(a_1, \ldots, a_M) \in S$, alors $(f_+(a_1), \ldots, f_+(a_M)) \in S$ et $(f_-(a_1), \ldots, f_-(a_M)) \in S$, ce qui permet de définir des probabilités de transition sur $S$ par

$$p((a_1, \ldots, a_M), (f_+(a_1), \ldots, f_+(a_M))) = \tfrac{1}{2},$$

$$p((a_1, \ldots, a_M), (f_-(a_1), \ldots, f_-(a_M))) = \tfrac{1}{2}.$$

Soit $C \subset S$ une classe récurrente, et soit $(X_n^1, \ldots, X_n^M)_{n \in \mathbf{Z}}$ une chaîne de Markov stationnaire sur $C$ pour les probabilités de transition $p$. Les arguments de la Proposition 5 montrent qu'il existe des jeux de pile ou face $(\xi_n^1)_{n \in \mathbf{Z}}, \ldots, (\xi_n^M)_{n \in \mathbf{Z}}$ ayant même image par $\Phi$ tels que presque sûrement, $X_n^k = (\xi_{n-d+1}^k, \ldots, \xi_n^k)$ pour tout $k \in [1 \ldots M]$ et pour tout $n \in \mathbf{Z}$. Ces jeux de pile ou face sont presque sûrement distincts deux-à-deux, puisque pour tout $n \in \mathbf{Z}$, les $d$-uplets $X_n^k = (\xi_{n-d+1}^k, \ldots, \xi_n^k)$ sont tous différents.

Notons $\mu$ la loi uniforme sur $\{-1, 1\}^{\mathbf{Z}}$. Pour tout $\alpha \in \{-1, 1\}^{\mathbf{Z}}$, notons $\nu(\alpha, \cdot)$ la loi de $\varepsilon$ sachant $\Phi(\varepsilon) = \alpha$. Soit $\zeta$ l'image commune par $\Phi$ des jeux de pile ou face $\xi^1, \ldots, \xi^M$. Alors pour tout $k \in [1 \ldots M]$ le couple $(\zeta, \xi^k)$ a même loi que $(\Phi(\varepsilon), \varepsilon)$. Comme la loi $\nu(\alpha, \cdot)$ est atomique pour $\mu$-presque tout $\alpha \in \{-1, 1\}^{\mathbf{Z}}$, on en déduit que presque sûrement, $\xi^1, \ldots, \xi^M$ sont des atomes de la loi $\nu(\zeta, \cdot)$. Donc $\nu(\alpha, \cdot)$ possède au moins $M$ atomes pour $\mu$-presque tout $\alpha$. $\quad\square$

A partir du Théorème 20, il est très facile d'obtenir un complément indépendant à $\mathcal{F}_n^{H\varepsilon}$ pour obtenir $\mathcal{F}_n^{\varepsilon}$:

**COROLLAIRE 22.** *Soit* $F : \{-1, 1\}^{-\mathbf{N}} \to \mathcal{C}$ *l'application fournie par la Proposition* 18. *Pour tout* $A \in \mathcal{C}$, *numérotons de* 1 *à* $M$ *les éléments de* $A$. *Pour tout* $n \in \mathbf{Z}$, *notons* $N_n$ *le numéro de* $(\varepsilon_{n-d+1}, \ldots, \varepsilon_n)$ *dans* $F((H\varepsilon)_{n]})$. *Alors presque sûrement,*

- $N_n$ *suit la loi uniforme sur* $\{1,\dots,M\}$;
- $N_n$ *est indépendante de* $\mathcal{F}_\infty^{H\varepsilon}$;
- $\mathcal{F}_n^\varepsilon = \mathcal{F}_n^{H\varepsilon} \vee \sigma(N_n)$.

REMARQUE. Il est facile de reconstituer la suite $\varepsilon_{n]}$ à partir de la suite $(H\varepsilon)_{n]}$ et de la variable aléatoire $N_n$. En effet, $(\varepsilon_{n-d+1},\dots,\varepsilon_n)$ est l'élément de numéro $N_n$ dans l'ensemble $F((H\varepsilon)_{n]})$. Par ailleurs, pour toute partie $A \in \mathcal{C}$, $f_+$ et $f_-$ induisent des bijections de $A$ dans leur image, puisque $A$, $f_+(A)$ et $f_-(A)$ sont de cardinal $M$. Connaissant $\varepsilon_{n-d+1},\dots,\varepsilon_n$ et $(H\varepsilon)_{n]}$, la relation $(\varepsilon_{n-d+1},\dots,\varepsilon_n) = f_{H_n\varepsilon_n}(\varepsilon_{n-d},\dots,\varepsilon_{n-1})$ détermine complètement l'élément $(\varepsilon_{n-d},\dots,\varepsilon_{n-1})$ dans l'ensemble $F((H\varepsilon)_{n-1]})$. On continue par récurrence décroissante.

**2.6. L'information perdue est un nombre fini de bits indépendants.** Nous avons vu dans la démonstration du Théorème 20 que la loi de $\varepsilon$ sachant $\Phi(\varepsilon)$ est (presque sûrement) uniforme sur un ensemble à $M$ éléments, où $M$ est le nombre maximum d'éléments de $E_d$ deux-à-deux non accordables. Nous allons voir que ce nombre $M$ est nécessairement une puissance de 2. Par conséquent, l'information perdue par la transformation $\Phi$ peut être représentée par un nombre fini de bits indépendants.

THÉORÈME 23. *Il existe un entier* $k \in [0\dots d]$ *tel que* $2^k$ *est le nombre maximum d'éléments de* $E_d$ *deux à deux non accordables, et tel que pour tout* $a \in E_d$, $2^{d-k}$ *est le nombre maximum d'éléments de* $E_d$ *simultanément accordables avec* $a$.

La démonstration du Théorème 23 s'appuie sur les deux propositions ci-dessous.

PROPOSITION 24. *Pour tout* $a \in E_d$, *soit* $N(a)$ *le nombre maximum d'éléments de* $E_d$ *simultanément accordables avec* $a$. *Alors* $N(a)$ *ne dépend pas de* $a$.

DÉMONSTRATION. Soit $a_0 \in E_d$ tel que $N(a_0)$ soit maximum. Soit $F$ un ensemble de $N(a_0)$ éléments simultanément accordables avec $a_0$. Comme chaque élément de $E_d$ a exactement deux antécédents par $f_+$ et $f_-$ réunis, on a donc $\mathrm{Card}(f_+^{-1}(F)) + \mathrm{Card}(f_-^{-1}(F)) = 2\,\mathrm{Card}(F)$. Mais $f_+^{-1}(F)$ et $f_-^{-1}(F)$ sont constituées chacune de points simultanément accordables et $\mathrm{Card}(F)$ est le nombre maximum d'éléments de $E_d$ simultanément accordables. On en déduit que $\mathrm{Card}(f_+^{-1}(F)) = \mathrm{Card}(f_-^{-1}(F)) = \mathrm{Card}(F)$, ce qui montre que $N(a) = N(a_0)$ pour les deux antécédents de $a_0$ par $f_+$ et $f_-$. Comme on peut atteindre $a_0$ par $d$ images directes par $f_+$ ou $f_-$ depuis n'importe quel point de $E_d$, on en déduit que $N(a)$ est indépendant de $a$. $\square$

Dans toute la suite, on notera $N$ le nombre maximum d'éléments simultanément accordables de $E_d$, si bien que $N = N(a)$ pour tout $a \in E_d$.

PROPOSITION 25.   *Soit $\alpha_1, \ldots, \alpha_l$ une suite finie de signes telle que la composée $f_{\alpha_l} \circ \cdots \circ f_{\alpha_1}$ accorde $N$ éléments en un élément $a_1$ et telle que l'ensemble $A = f_{\alpha_l} \circ \cdots \circ f_{\alpha_1}(E_d)$ soit formé de $M$ éléments deux-à-deux non accordables. Alors les images réciproques des éléments de $A$ par $f_{\alpha_l} \circ \cdots \circ f_{\alpha_1}$ sont toutes de cardinal $N$.*

L'existence d'une suite finie de signes vérifiant les conditions de la Proposition 25 découle des arguments de la partie 2.3: on commence par accorder simultanément $N$ points, puis on compose au hasard par $f_+$ et $f_-$. Le théorème découle immédiatement de la Proposition 25. En effet, comme les images réciproques des éléments de $A$ par $f_{\alpha_l} \circ \cdots \circ f_{\alpha_1}$ forment une partition de $E_d$, on en déduit immédiatement la relation $MN = 2^d$.

DÉMONSTRATION.   Soient $a_1, \ldots, a_M$ les éléments de $A$. Pour tout $k \in [1 \ldots M]$, notons $F_k = (f_{\alpha_l} \circ \cdots \circ f_{\alpha_1})^{-1}(a_k)$ et $N_k = \mathrm{Card}(F_k)$. Alors $N_k \leq N$ puisque les éléments de $F_k$ sont simultanément accordables, ce qui entraîne $2^d = N_1 + \cdots + N_M \leq MN$.

Pour montrer que $N_k = N$ pour tout $k \in [1 \ldots M]$, il suffit de montrer que $2^d = MN$. Pour cela, on note $X_n = (\varepsilon_{n-d+1}, \ldots, \varepsilon_n)$ pour tout $n \in \mathbf{Z}$ et on remarque que

$$X_l = (f_{H_l \varepsilon_l} \circ \cdots \circ f_{H_1 \varepsilon_1})(X_0),$$

d'où, par indépendance de $X_0$ et $(H_1 \varepsilon_1, \ldots, H_l \varepsilon_l)$,

$$P[X_l = a_1 | (H_1 \varepsilon_1, \ldots, H_l \varepsilon_l) = (\alpha_1, \ldots, \alpha_l)]$$

$$= P[X_0 \in F_1 | (H_1 \varepsilon_1, \ldots, H_l \varepsilon_l) = (\alpha_1, \ldots, \alpha_l)] = P[X_0 \in F_1] = \frac{N}{2^d}.$$

Donc $P[(\varepsilon_{l-d+1}, \ldots, \varepsilon_l) = a_1 | H \varepsilon] \geq \frac{N}{2^d}$ avec probabilité strictement positive. Mais d'après le Théorème 20 cette probabilité conditionnelle est (presque sûrement) à valeurs dans $\{0, \frac{1}{M}\}$. Donc $\frac{1}{M} \geq \frac{N}{2^d}$, d'où $2^d \geq MN$, ce qui donne l'égalité voulue.   □

Maintenant que la relation $MN = 2^d$ est établie, on peut a posteriori affaiblir les hypothèses de la Proposition 25.

COROLLAIRE 26.   *Soit $f_{\alpha_l} \circ \cdots \circ f_{\alpha_1}$ une composée de $f_+$ et $f_-$ réalisant un accord maximal, c'est-à-dire telle que l'ensemble $A = f_{\alpha_l} \circ \cdots \circ f_{\alpha_1}(E_d)$ soit formé de $M$ éléments deux-à-deux non accordables. Alors les images réciproques des éléments de $A$ par $f_{\alpha_l} \circ \cdots \circ f_{\alpha_1}$ sont toutes de cardinal $N$.*

Démonstration. Il suffit de reprendre le début de la démonstration de la Proposition 25 en remarquant que les inégalités $N_k \leq N$ et l'égalité $N_1 + \cdots + N_M = 2^d = MN$ entraînent $N_k = N$ pour tout $k$. □

Corollaire 27. *Il existe une partition de $E_d$ en $M$ parties disjointes formées chacune de $N$ éléments simultanément accordables.*

Dans la démonstration de l'égalité $MN = 2^d$ que nous avons donnée, la fin (démonstration de l'inégalité $MN \leq 2^d$) utilise les arguments probabilistes du Théorème 20. Nous proposons ci-dessous une démonstration combinatoire directe et constructive du Corollaire 27, et par là-même du Théorème 23. La démonstration repose sur les deux lemmes suivants:

Lemme 28. *Soit $F$ une partie de $E_d$ formée de $N$ éléments simultanément accordables. Alors l'image réciproque de $F$ par $f_+$ et $f_-$ est formée de $N$ éléments simultanément accordables.*

Démonstration. Il suffit de reprendre le début la démonstration de Proposition 24 en remarquant que les inégalités $\mathrm{Card}(f_+^{-1}(F)) \leq N$ et $\mathrm{Card}(f_-^{-1}(F)) \leq N$, et l'égalité $\mathrm{Card}(f_+^{-1}(F)) + \mathrm{Card}(f_-^{-1}(F)) = 2\,\mathrm{Card}(F)$ entraînent $\mathrm{Card}(f_+^{-1}(F)) = \mathrm{Card}(f_-^{-1}(F)) = N$. □

Lemme 29. *Soient $F_1, \ldots, F_k$ des parties disjointes de $E_d$ formée chacune de $N$ éléments simultanément accordables, et $g$ une composée de $f_+$ et $f_-$ constante sur chacune de ces parties. Notons $c_1, \ldots, c_k$ la valeur prise par $g$ sur les parties $F_1, \ldots, F_k$. Alors:*

- *Les éléments $c_1, \ldots, c_k$ sont deux-à-deux non accordables, et donc $k \leq M$.*
- *Si $F_1 \cup \cdots \cup F_k \neq E$ (en particulier si $k < M$), on peut construire $k+1$ parties disjointes de $E_d$ formée chacune de $N$ éléments simultanément accordables, et une composée de $f_+$ et $f_-$ constante sur chacune de ces parties.*

Démonstration. Le premier point est immédiat: s'il existait deux indices $i < j$ tels que les points $c_i$ et $c_j$ soient accordables par une composée $h$, alors la composée $h \circ g$ accorderait les $2M$ éléments de $F_i \cup F_j$.

Montrons le second point: soit donc $a \in E_d \setminus (F_1 \cup \cdots \cup F_k)$ et $c_{k+1} = g(a)$. Alors $c_{k+1}$ est différent de $c_1, \ldots, c_k$ (sinon, on obtiendrait $N+1$ éléments simultanément accordables en ajoutant $a$ à l'une des parties $F_1, \ldots, F_k$). Soit $h$ une composée de $f_+$ et $f_-$ envoyant $c_1$ sur $a$. Alors les images réciproques de $c_1, \ldots, c_k, c_{k+1}$ par $g \circ h \circ g$ sont formées chacune de $N$ éléments simultanément accordables. □

Une application répétée du Lemme 29 fournit une partition de $E_d$ en $M$ parties disjointes formées chacune de $N$ éléments simultanément accordables, ce qui démontre le Corollaire 27 et le Théorème 23.

2.7. *Majoration du nombre de bits perdus et condition suffisante de conservation de l'information.* A l'aide du Lemme 21 et du Théorème 23, on voit immédiatement que la transformation $\Phi : (x_n)_{n \in \mathbf{Z}} \mapsto (\phi(x_{n-d}, \ldots, x_{n-1})x_n)_{n \in \mathbf{Z}}$ fait perdre au plus $d$ bits d'informations. On peut améliorer cette borne en faisant quelques hypothèses simples sur la loi du processus stationnaire $(H_n)_{n \in \mathbf{Z}}$ défini par $H_n = \phi(\varepsilon_{n-d}, \ldots, \varepsilon_{n-1})$.

THÉORÈME 30. *Soit $p \in [0 \ldots d]$. S'il existe $(h_1, \ldots, h_d) \in \{-1, 1\}^d$ tel que $P[(H_1, \ldots, H_d) = (h_1, \ldots, h_d)] > \frac{1}{2^{p+1}}$, la transformation $\Phi$ perd au plus $p$ bits d'information. En particulier, si $P[(H_1, \ldots, H_d) = (1, \ldots, 1)] > \frac{1}{2}$ ou si $P[(H_1, \ldots, H_d) = (-1, \ldots, -1)] > \frac{1}{2}$, la transformation $\Phi$ conserve l'information.*

DÉMONSTRATION. Supposons que $P[(H_1, \ldots, H_d) = (h_1, \ldots, h_d)] > \frac{1}{2^{p+1}}$. On sait d'après le Théorème 20 que la loi de $(\varepsilon_1, \ldots, \varepsilon_d)$ sachant $\Phi(\varepsilon)$ est la loi uniforme sur un ensemble à $M$ éléments. Par conséquent,

$$P[(H_1, \ldots, H_d) = (h_1, \ldots, h_d)]$$
$$= P[(\varepsilon_1, \ldots, \varepsilon_d) = (h_1 H_1 \varepsilon_1, \ldots, h_d H_d \varepsilon_d)] \leq \frac{1}{M},$$

puisque $(h_1 H_1 \varepsilon_1, \ldots, h_d H_d \varepsilon_d)$ est une fonction de $\Phi(\varepsilon)$. On a donc $\frac{1}{2^{p+1}} < \frac{1}{M}$. Comme $M$ est une puissance de 2, cela entraîne $M \leq 2^p$. $\square$

REMARQUE. Comme le $d$-uplet $(H_1, \ldots, H_d)$ peut prendre au plus $2^d$ valeurs, on peut toujours appliquer le Théorème 30 avec $p = d$. Nous verrons au paragraphe 3.2 une condition nécessaire et suffisante pour que la transformation $\Phi$ perde exactement $d$ bits d'information.

Le corollaire ci-dessous montre que si la variable $H_1$ est trop proche d'une variable aléatoire déterministe, la transformation $\Phi$ conserve l'information.

COROLLAIRE 31. *Si $P[H_1 \neq 1] < \frac{2^{p+1}-1}{2^{p+1}d}$, la transformation $\Phi$ perd au plus $p$ bits d'information. En particulier, si $P[H_1 \neq 1] < \frac{1}{2d}$, la transformation $\Phi$ conserve l'information.*

DÉMONSTRATION. Il suffit d'écrire, par équidistribution des variables $H_1, \ldots, H_d$,

$$P[(H_1, \ldots, H_d) \neq (1, \ldots, 1)] \leq d P[H_1 \neq 1] < \frac{2^{p+1}-1}{2^{p+1}} = 1 - \frac{1}{2^{p+1}},$$

et d'appliquer le théorème précédent. $\square$

**3. Transformations homogènes de longueur finie: exemples.** Dans cette partie, nous donnons des conditions suffisantes de perte ou de conservation de l'information, ainsi que de nombreux exemples illustrant les résultats de la deuxième partie.

### 3.1. *Longueur des transformations.*

DÉFINITION 32. On appelle transformation homogène de longueur finie toute application $\Phi : \{-1; 1\}^{\mathbf{Z}} \to \{-1; 1\}^{\mathbf{Z}}$ de la forme

$$\Phi(x) = (\phi(x_{n-d}, \dots, x_{n-1})x_n)_{n \in \mathbf{Z}} \qquad \text{pour tout } x = (x_n)_{n \in \mathbf{Z}} \in \{-1; 1\}^{\mathbf{Z}}$$

avec $d \in \mathbf{N}$ et $\phi : \{-1; 1\}^d \to \{-1; 1\}$. Dans ce cas, le plus petit entier naturel $d$ tel qu'on puisse écrire $\Phi$ sous cette forme s'appelle longueur de $\Phi$.

Exemples: les seules transformations de longueur 0 sont Id et $-$Id. Les seules transformations de longueur 1 sont $\Phi_1 : (x_n)_{n \in \mathbf{Z}} \to (x_n x_{n-1})_{n \in \mathbf{Z}}$ et $-\Phi_1$.

PROPOSITION 33. *Si $\Phi$ et $\Psi$ sont des transformations homogènes de longueur $p$ et $q$, alors $\Psi \circ \Phi$ est une transformation homogène de longueur $p + q$.*

DÉMONSTRATION. Par hypothèse, il existe des applications $\phi : \{-1; 1\}^p \to \{-1; 1\}$ et $\psi : \{-1; 1\}^q \to \{-1; 1\}$ telles que pour tout $x = (x_n)_{n \in \mathbf{Z}} \in \{-1; 1\}^{\mathbf{Z}}$,

$$\Phi(x) = (\phi(x_{n-p}, \dots, x_{n-1})x_n)_{n \in \mathbf{Z}},$$
$$\Psi(x) = (\psi(x_{n-q}, \dots, x_{n-1})x_n)_{n \in \mathbf{Z}}.$$

Soit $x = (x_n)_{n \in \mathbf{Z}} \in \{-1; 1\}^{\mathbf{Z}}$, $y = (y_n)_{n \in \mathbf{Z}} = \Phi(x)$ et $z = (z_n)_{n \in \mathbf{Z}} = \Psi(y)$. Alors pour tout $n \in \mathbf{Z}$,

$$\begin{aligned}
z_n &= \psi(y_{n-q}, \dots, y_{n-1})y_n \\
&= \psi(\phi(x_{n-p-q}, \dots, x_{n-1-q})x_{n-q}, \dots, \phi(x_{n-p-1}, \dots, x_{n-2})x_{n-1}) \\
&\quad \times \phi(x_{n-p}, \dots, x_{n-1})x_n \\
&= \theta(x_{n-p-q}, \dots, x_{n-1})x_n,
\end{aligned}$$

en notant pour tout $(t_1, \dots, t_{p+q}) \in \{-1, 1\}^{p+q}$:

$$\begin{aligned}
&\theta(t_1, \dots, t_{p+q}) \\
&\quad = \psi(\phi(t_1, \dots, t_p)t_{p+1}, \dots, \phi(t_q, \dots, t_{p+q-1})t_{p+q})\phi(t_{1+q}, \dots, t_{p+q}).
\end{aligned}$$

Donc la transformation $\Psi \circ \Phi$ est homogène de longueur $\leq p + q$. Pour montrer qu'elle est vraiment de longueur $p+q$, il reste à vérifier que $\theta(t_1, \dots, t_{p+q})$ dépend effectivement de $t_1$.

Par hypothèse, les applications $\phi$ et $\psi$ dépendent effectivement de la première variable, autrement dit, il existe $\alpha_2, \ldots, \alpha_p, \beta_2, \ldots, \beta_q \in \{-1, 1\}$ tels que

$$\phi(-1, \alpha_2, \ldots, \alpha_p) \neq \phi(1, \alpha_2, \ldots, \alpha_p),$$
$$\psi(-1, \beta_2, \ldots, \beta_q) \neq \psi(1, \beta_2, \ldots, \beta_q).$$

Définissons $\gamma_2, \ldots, \gamma_{p+q}$ en posant $\gamma_k = \alpha_k$ pour tout $k \in [2 \ldots p]$, $\gamma_{p+1} = \phi(1, \alpha_2, \ldots, \alpha_p)$, et $\gamma_{p+k} = \phi(\gamma_k, \ldots, \gamma_{p+k-1})\beta_k$ pour tout $k \in [2 \ldots q]$. Soit $c = \phi(\gamma_{1+q}, \ldots, \gamma_{p+q})$. Alors

$$\begin{aligned}
\theta(-1, \gamma_2, \ldots, \gamma_{p+q}) &= c\psi(-1, \beta_2, \ldots, \beta_q) \\
&\neq c\psi(1, \beta_2, \ldots, \beta_q) \\
&= \theta(1, \gamma_2, \ldots, \gamma_{p+q})
\end{aligned}$$

ce qui montre que $\theta$ dépend effectivement de la première variable. $\square$

Nous avons vu dans la démonstration directe et constructive de l'implication $(3) \Rightarrow (1)$ du Théorème 11 que si $\Phi$ est une transformation homogène de longueur finie qui conserve l'information, elle possède un inverse $\Psi$ qui est également une transformation homogène. Mais par additivité des longueurs, la transformation $\Psi$ ne peut être de longueur finie, à moins que $\Phi$ soit de longueur 0.

COROLLAIRE 34. *L'inverse d'une transformation homogène de longueur finie qui conserve l'information, autre que $\pm \mathrm{Id}$, est une transformation homogène de longueur infinie.*

3.2. *Conditions suffisantes de perte d'information.* Nous allons voir trois exemples simples de transformations faisant perdre de l'information: les transformations paires, les transformations "dichotomiques" et les transformations pour lesquelles $f_+$ et $f_-$ sont bijectives.

Dans la suite, on considère une transformation $\Phi : \{-1; 1\}^{\mathbf{Z}} \to \{-1; 1\}^{\mathbf{Z}}$ de longueur $d$, associée à $\phi : \{-1; 1\}^d \to \{-1; 1\}$. On note $\Phi_1 : (x_n)_{n \in \mathbf{Z}} \to (x_n x_{n-1})_{n \in \mathbf{Z}}$. La transformation $\Phi_1$ est de longueur 1 et associée à $\phi_1 : x_1 \mapsto x_1$. C'est l'exemple le plus simple de transformation faisant perdre de l'information.

Nous allons voir que les transformations paires ou "dichotomiques" perdent de l'information simplement parce qu'elles s'écrivent comme composées (dans un sens ou dans l'autre) de $\Phi_1$ et d'une autre transformation.

Les transformations pour lesquelles $f_+$ et $f_-$ sont bijectives, quant à elles, sont celles qui perdent le maximum d'information, soit un nombre de bits égal à leur longueur. Par la suite, nous les appellerons "transformations à perte maximum."

*Les transformations paires.*

PROPOSITION 35.    *Il y a équivalence entre:*

(1)  $\Phi$ *est paire;*
(2)  $\phi$ *est impaire;*
(3)  *il existe une application* $\psi \colon \{-1;1\}^{d-1} \to \{-1;1\}$ *telle que pour tout* $(x_1, \ldots, x_d) \in \{-1;1\}^d$, $\phi(x_1, \ldots, x_d) = x_1 \psi(x_1 x_2, \ldots, x_{d-1} x_d);$
(4)  *il existe une transformation* $\Psi$ *de longueur* $d-1$ *telle que* $\Phi = \Psi \circ \Phi_1$.

DÉMONSTRATION.    Les implications $(1) \Rightarrow (2) \Rightarrow (3) \Rightarrow (4) \Rightarrow (1)$ sont immédiates.    □

Exemple: prenons $d = 3$ et $\phi(x_1, x_2, x_3) = \operatorname{sgn}(x_1 + x_2 + x_3)$. Alors par imparité de $\phi$, pour tout $(x_1, x_2, x_3) \in \{-1;1\}^3$,

$$\phi(x_1, x_2, x_3) = x_1 \operatorname{sgn}(1 + x_1 x_2 + x_1 x_3)$$
$$= x_1 \max(x_1 x_2, x_1 x_3)$$
$$= x_1 \max(x_1 x_2, -x_2 x_3).$$

Donc $\Phi = \Psi \circ \Phi_1$, où $\Psi$ est la transformation de longueur 2 associée à $\psi \colon (y_1, y_2) \mapsto \max(y_1, -y_2)$.

*Les transformations dichotomiques.*

PROPOSITION 36.    *Soit* $\gamma \in \{-1, 1\}$. *Il y a équivalence entre:*

(1)  *Il existe une partition de* $\{-1, 1\}^d$ *en deux parties* $A$ *et* $B$ *telle que* $f_\gamma(A) \subset A$, $f_\gamma(B) \subset B$, $f_{-\gamma}(A) \subset B$ *et* $f_{-\gamma}(B) \subset A;$
(2)  *il existe une application* $\psi \colon \{-1;1\}^{d-1} \to \{-1;1\}$ *telle que pour tout* $(x_1, \ldots, x_d) \in \{-1;1\}^d$, $\phi(x_1, \ldots, x_d) = \gamma x_d \psi(x_1, \ldots, x_{d-1}) \psi(x_2, \ldots, x_d);$
(3)  *il existe une transformation* $\Psi$ *de longueur* $d-1$ *telle que* $\Phi = \gamma(\Phi_1 \circ \Psi)$.

*Lorsque ces conditions sont vérifiées,* $A$ *et* $B$ *sont les graphes de* $\psi$ *et* $-\psi$ *(ou inversement) et* $\gamma = \phi(-1, \ldots, -1) = \phi(1, \ldots, 1)$. *De plus, si* $A$ *est le graphe de* $\psi$, *chaque variable aléatoire* $X_n = \mathbf{I}_A(\varepsilon_{n-d+1}, \ldots, \varepsilon_n) - \mathbf{I}_B(\varepsilon_{n-d+1}, \ldots, \varepsilon_n) = \varepsilon_n \psi(\varepsilon_{n-d+1}, \ldots, \varepsilon_{n-1})$ *est indépendante de* $\Phi(\varepsilon)$.

DÉMONSTRATION.    On démontre les implications $(1) \Rightarrow (2) \Rightarrow (3) \Rightarrow (1)$, la seule implication non immédiate étant $(1) \Rightarrow (2)$. Supposons que la condition (1) soit vérifiée. Alors pour tout $(x_1, \ldots, x_{d-1}) \in \{-1;1\}^{d-1}$, l'un des $d$-uplets $(x_1, \ldots, x_{d-1}, 1)$ et $(x_1, \ldots, x_{d-1}, -1)$ appartient à $A$ et l'autre à $B$, puisque ce sont les images par $f_+$ et $f_-$ (ou inversement) du $d$-uplet $(1, x_1, \ldots, x_{d-1})$. Par conséquent, $A$ est le graphe d'une application $\psi \colon \{-1;1\}^{d-1} \to \{-1;1\}$ et $B$ est le graphe de $-\psi$. Mais pour tout $x = (x_1, \ldots, x_d) \in \{-1;1\}^d$,

$$(x_2, \ldots, x_d, \gamma\phi(x_1, \ldots, x_d)) = f_\gamma(x_1, \ldots, x_d)$$

et $f_\gamma(x) \in A \Leftrightarrow x \in A$. On en déduit que

$$\gamma\phi(x_1, \ldots, x_d) = \psi(x_2, \ldots, x_d) \quad \Longleftrightarrow \quad x_d = \psi(x_1, \ldots, x_{d-1}).$$

Comme toutes les quantités impliquées valent $-1$ ou $1$, on a donc

$$\gamma\phi(x_1, \ldots, x_d)\psi(x_2, \ldots, x_d) = x_d\psi(x_1, \ldots, x_{d-1}),$$

ce qui montre la condition (2). $\quad \square$

REMARQUE. La condition (3) semble beaucoup plus simple que (1) et explique en particulier pourquoi une transformation dichotomique perd de l'information. Mais d'un point de vue pratique, il est souvent bien plus facile de trouver la partition de $\{-1,1\}^d$ en $A$ et $B$ que de deviner la fonction $\psi$, ou ce qui revient au même la transformation $\Psi$.

Exemple: prenons $d = 3$ et $\phi(x_1, x_2, x_3) = \min(x_1 x_3, x_2)$. On voit immédiatement que la transformation $\Phi$ est dichotomique puisque $A = \{(-, +, +), (+, +, -), (+, -, -), (-, -, -)\}$ et $B = \{(-, -, +), (-, +, -), (+, -, +), (+, +, +)\}$ sont stables par $f_+$, tandis que $f_-$ envoie $A$ dans $B$ et $B$ dans $A$.

La partie $A$ est le graphe de $\psi \colon (x_1, x_2) \mapsto \min(-x_1, x_2)$. Pour tout $(x_1, x_2, x_3) \in \{-1;1\}^3$, on a en effet $\phi(x_1, x_2, x_3) = x_3\min(-x_1, x_2)\min(-x_2, x_3)$. La vérification de cette égalité est immédiate (en discutant sur la valeur de $x_2$), mais il n'était pas évident de deviner cette décomposition a priori.

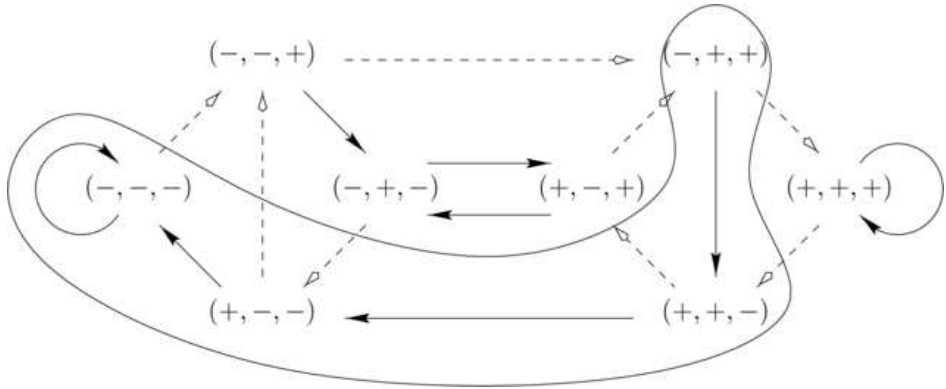

FIG. 3. *Les flèches pleines indiquent l'image par $f_+$, les flèches creuses l'image par $f_-$.*

*Les transformations à perte maximum.*

PROPOSITION 37. *Il y a équivalence entre:*

(1) *La transformation $\Phi$ fait perdre $d$ bits d'information;*
(2) *Les éléments de $E_d$ sont deux-à-deux non accordables;*
(3) *Les applications $f_+$ et $f_-$ sont bijectives;*
(4) *il existe une application $\psi\colon\{-1;1\}^{d-1}\to\{-1;1\}$ telle que pour tout $(x_1,\ldots,x_d)\in\{-1;1\}^d$, $\phi(x_1,\ldots,x_d)=x_1\psi(x_2,\ldots,x_d)$;*

*Lorsque ces conditions sont vérifiées, chaque $d$-uplet $(\varepsilon_{n-d+1},\ldots,\varepsilon_n)$ est indépendant de $\Phi(\varepsilon)$.*

DÉMONSTRATION. L'implication $(1)\Rightarrow(2)$ découle des résultats du Paragraphe 2.6. Les implications $(2)\Rightarrow(3)\Rightarrow(4)$ sont immédiates.

Montrons l'implication $(4)\Rightarrow(1)$ et la dernière affirmation: supposons que $\phi$ soit de la forme $\phi(x_1,\ldots,x_d)=x_1\psi(x_2,\ldots,x_d)$. Soit $n\in\mathbf{Z}$. Alors quelles que soient les parties finies $I$ de $\mathbf{N}$ et $J$ non vide de $[1\ldots d]$,

$$\mathbf{E}\left[\prod_{i\in I}(H_{n-i}\varepsilon_{n-i})\prod_{j\in J}\varepsilon_{n-j}\right]=0,$$

par indépendance des variables $\varepsilon_k$. En effet, le résultat est évident si $I$ est vide. Et si $I$ est non vide et si $m=\max I$, la variable dont on calcule l'espérance s'écrit comme le produit de $\varepsilon_{n-m-d}$ par une fonction des variables $\varepsilon_k$ pour $n-m-d+1\le k\le n$. Par conséquent, le $d$-uplet $(\varepsilon_{n-d+1},\ldots,\varepsilon_n)$ est indépendant de $\Phi(\varepsilon)_{n]}$, et même de $\Phi(\varepsilon)$ par indépendance de $(\Phi(\varepsilon)_k)_{k\ge n+1}$ et de $\mathcal{F}_n^\varepsilon$. Il suffit alors d'utiliser les résultats du Paragraphe 2.4.  □

On peut se demander si toute transformation de longueur finie perdant de l'information est la composée d'une transformation à perte maximum et d'autres transformations. Il n'en n'est rien, et nous verrons au paragraphe 3.4 un exemple de transformation indécomposable perdant de l'information mais qui n'est pas une transformation à perte maximum.

3.3. *Un peu de dénombrement.* Soit $d\in\mathbf{N}^*$. Le nombre de transformations de longueur $\le d$ est le nombre d'applications de $\phi\colon\{-1;1\}^d\to\{-1;1\}$, soit $2^{2^d}$. Parmi ces transformations, comptons le nombre de transformations simples faisant perdre de l'information.

Il y a exactement $2^{2^{d-1}}$ transformations paires: autant que d'applications $\phi$ impaires de $\{-1;1\}^d$ dans $\{-1;1\}$.

Il y a exactement $2^{2^{d-1}}$ transformations dichotomiques: autant que d'applications $\psi$ de $\{-1;1\}^{d-1}$ dans $\{-1;1\}$. En effet, notons $T$ le morphisme de groupes qui à une application $\psi$ de $\{-1;1\}^{d-1}$ dans $\{-1;1\}$ associe l'application $T(\psi)\colon\{-1;1\}^d \to \{-1;1\}$ définie par $T(\psi)(x_1,\dots,x_d) = \psi(x_1,\dots,x_{d-1})\psi(x_2, \dots,x_d)$. On vérifie facilement que le noyau de $T$ est formé des applications constantes égales à 1 ou $-1$. Par conséquent, à deux applications opposées de $\{-1;1\}^{d-1}$ dans $\{-1;1\}$ correspondent deux transformations dichotomiques, obtenues en prenant $\gamma = 1$ ou $\gamma = -1$.

Il y a exactement $2^{2^{d-1}}$ transformations de longueur $d$ à perte maximum: autant que d'applications $\psi$ de $\{-1;1\}^{d-1}$ dans $\{-1;1\}$.

Certaines transformations peuvent appartenir à plusieurs catégories à la fois. Par exemple, la transformation $\Phi_1 \circ \Phi_1$, de longueur 2, est à la fois paire, dichotomique et à perte maximum. Donnons les cardinaux des différentes intersections pour $d \geq 3$ (voir figure 4).

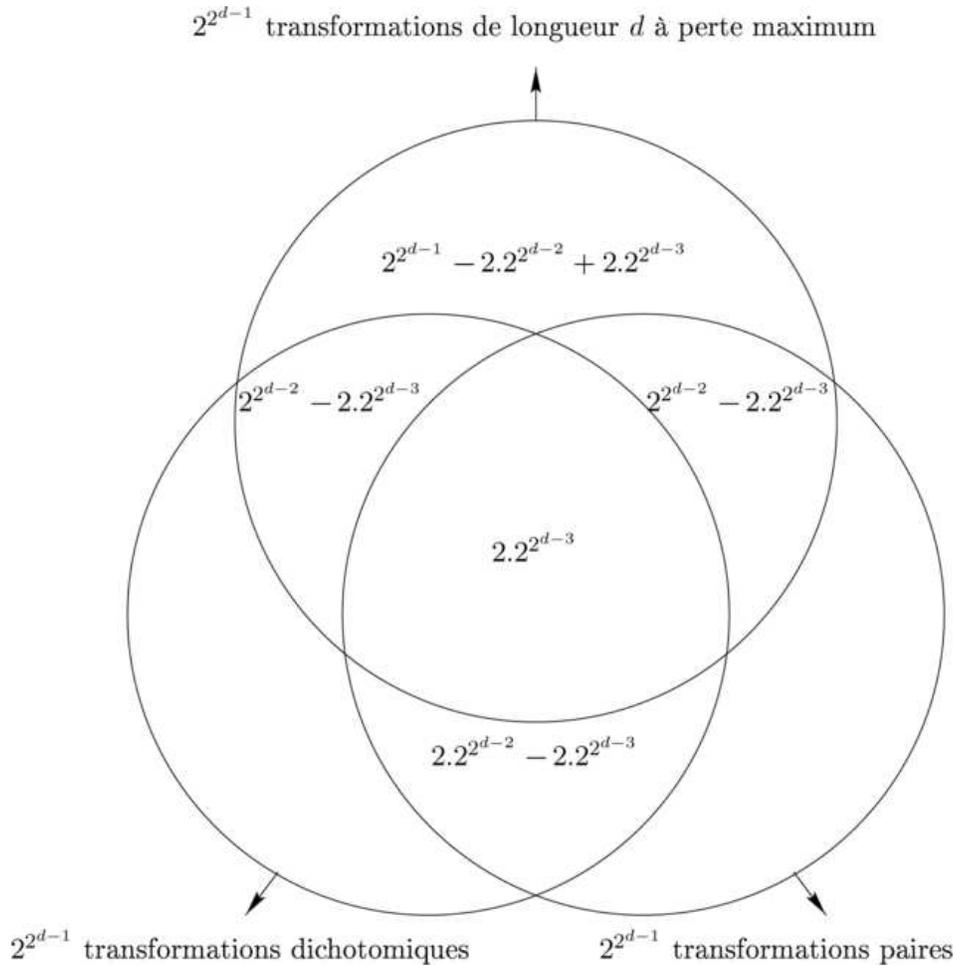

Fig. 4.

Il y a exactement $2^{2^{d-2}}$ transformations paires de longueur $d$ à perte maximum: autant que d'applications $\phi$ de la forme $\phi(x_1,\ldots,x_d) = x_1\psi(x_2,\ldots,x_d)$ avec $\psi\colon\{-1;1\}^{d-1}\to\{-1;1\}$ paire (pour que $\phi$ soit impaire).

Il y a exactement $2.2^{2^{d-2}}$ transformations paires et dichotomiques: en effet, une application $\phi$ de la forme $\phi(x_1,\ldots,x_d) = \gamma x_d\psi(x_1,\ldots,x_{d-1})\psi(x_2,\ldots,x_d)$ avec $\gamma\in\{-1,1\}$ et $\psi\colon\{-1;1\}^{d-1}\to\{-1;1\}$ est paire si et seulement si $\psi$ est paire ou impaire: si on note $\widetilde{\psi}(x_1,\ldots,x_{d-1}) = \psi(-x_1,\ldots,-x_{d-1})$, alors $T(\widetilde{\psi}) = T(\psi)$, équivaut à $\widetilde{\psi}=\psi$ ou $\widetilde{\psi}=-\psi$.

Il y a exactement $2^{2^{d-2}}$ transformations dichotomiques de longueur $d$ à perte maximum: en effet, si $\phi$ est de la forme $\phi(x_1,\ldots,x_d) = \gamma x_d\psi(x_1,\ldots,x_{d-1}) \times \psi(x_2,\ldots,x_d)$ avec $\gamma\in\{-1,1\}$ et $\psi\colon\{-1;1\}^{d-1}\to\{-1;1\}$, alors $x_1\phi(x_1,\ldots,x_d)$ ne dépend pas de $x_1$ si et seulement si $x_1\psi(x_1,\ldots,x_{d-1})$ ne dépend pas de $x_1$, c'est-à-dire si $\psi$ est de la forme $\psi(x_1,\ldots,x_{d-1}) = x_1\theta(x_2,\ldots,x_{d-1})$ avec $\theta\colon\{-1;1\}^{d-2}\to\{-1;1\}$.

Parmi ces transformations, les transformations paires sont celles pour lesquelles l'application $\theta$ est paire ou impaire. Si $d\geq 3$, il y en a $2.2^{2^{d-3}}$. Le nombre de transformations qui sont paires, dichotomiques ou de longueur $d$ à perte maximum est donc $3.2^{2^{d-1}} - 4.2^{2^{d-2}} + 2.2^{2^{d-3}}$. Le nombre de transformations de longueur $d$ à perte maximum qui ne sont ni paires ni impaires est $2^{2^{d-1}} - 2.2^{2^{d-2}} + 2.2^{2^{d-3}}$.

3.4. *Etude pour les petites valeurs de $d$.*

*Le cas où $d=2$.*  Il est facile de passer en revue les 16 transformations de longueur $\leq 2$, ce qui donne les résultats suivants en fonction de la forme de l'application $\phi\colon\{-1;1\}^2\to\{-1;1\}$. Dans la liste ci-dessous $\alpha,\beta,\gamma$ sont des constantes valant $-1$ ou 1.

Si $\phi\colon(x_1,x_2)\mapsto\alpha$, alors $\Phi=\alpha\,\mathrm{Id}$ est de longueur 0 et conserve l'information.

Si $\phi\colon(x_1,x_2)\mapsto\alpha x_2$, alors $\Phi=\alpha\Phi_1$ est de longueur 1 et perd un bit d'information.

Si $\phi\colon(x_1,x_2)\mapsto\alpha x_1$, alors $\Phi=\alpha\Phi_1\circ\Phi_1$ est de longueur 2 et perd deux bits d'information.

Si $\phi\colon(x_1,x_2)\mapsto\alpha x_1 x_2$, alors $\Phi$ est de longueur 2 et perd deux bits d'information.

Si $\phi\colon(x_1,x_2)\mapsto\alpha\max(\beta x_1,\gamma x_2)$, alors $\Phi$ est de longueur 2 et conserve l'information.

REMARQUE.  Cette liste montre qu'une transformation $\Phi$ de longueur $\leq 2$ perd de l'information si et seulement si l'application $\phi\colon\{-1;1\}^2\to\{-1;1\}$ est "centrée," c'est-à-dire qu'elle prend autant de fois la valeur $-1$ que la valeur 1. Mais cette propriété ne subsiste pas pour des longueurs supérieures à 2, comme le montre l'exemple de la transformation dichotomique donnée au paragraphe précédent. Dans cet exemple, l'application $\phi \colon \{-1;1\}^3 \to \{-1;1\}$ prend 6 fois la valeur $-1$ et 2 fois la valeur 1.

*Le cas où $d = 3$.* L'algorithme de Roland Bacher appliqué aux 256 applications de $\{-1;1\}^3$ dans $\{-1;1\}$ fournit la liste des 38 transformations de longueurs $\leq 3$ qui perdent de l'information. Le dénombrement effectué au paragraphe précédent montre que ces 38 transformations sont toutes paires, dichotomiques ou à perte maximum, en incluant les deux transformations de longueur 2 à perte maximum $\Phi_2 \colon (x_n)_{n \in \mathbf{Z}} \to (x_n x_{n-1} x_{n-2})_{n \in \mathbf{Z}}$ et $-\Phi_2$, qui ne sont ni paires, ni dichotomiques. Parmi ces 38 transformations, seules quatre correspondent à des applications $\phi$ "non centrées" : la transformation dichotomique de l'exemple du Paragraphe 3.2 et ses trois composées avec $-\,\mathrm{Id}$ (à gauche, à droite et des deux côtés).

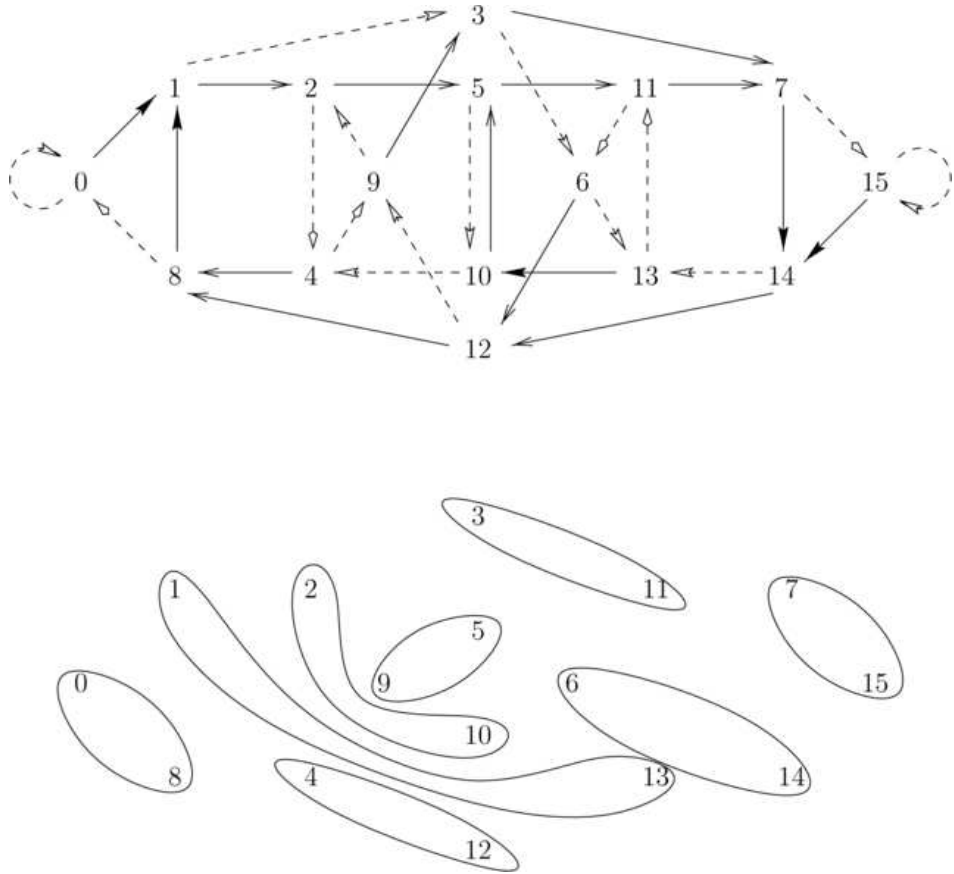

Fig. 5. *Graphe de $f_+$ et $f_-$, et relation d'accordabilité.*

*Le cas où $d = 4$.*   L'algorithme de Roland Bacher appliqué aux 65 536 applications de $\{-1; 1\}^3$ dans $\{-1; 1\}$ fournit la liste des 782 transformations de longueur $\leq 4$ qui perdent de l'information. Le dénombrement effectué au paragraphe précédent montre qu'il y a 712 transformations paires, dichotomiques ou de longueur 4 à perte maximum.

On peut ajouter 14 transformations à perte maximum de longueur $\leq$ 3 ni paires, ni dichotomiques: 2 de longueur 2 et 12 de longueur 3. On peut encore ajouter les composées de $\Phi_2$ et d'une autre transformation de longueur 2 qui n'ont pas encore été comptabilisées... Mais cela ne suffit pas à atteindre 782. Il y a donc des exemples "compliqués" de transformations qui perdent de l'information. Nous donnons ici un exemple de transformation indécomposable de longueur 4 qui perd 3 bits d'information.

Pour la clarté de la figure 5, on numérote de 0 à 15 les éléments de $\{-1; 1\}^4$ en donnant au quadruplet $(\alpha, \beta, \gamma, \delta)$ le numéro $8\frac{1+\alpha}{2} + 4\frac{1+\beta}{2} + 2\frac{1+\gamma}{2} + \frac{1+\delta}{2}$. Dans le graphe de $f_+$ et $f_-$, les flèches issues du quadruplet numéro $x$ sont alors dirigées vers les quadruplets numéro $2x$ et $2x + 1$ (modulo 16). L'application $\phi : \{-1; 1\}^4 \to \{-1; 1\}$ prend la valeur 1 aux quadruplets numéro 0, 2, 3, 5, 8, 9, 10, 11, et la valeur $-1$ sur les huit autres quadruplets.

On vérifie facilement que la relation d'accordabilité est une relation d'équivalence comportant 8 classes à 2 éléments chacune, si bien que la transformation $\Phi$ perd 3 bits d'information.

Comme $\phi(-1, -1, -1, 1) \neq \phi(1, -1, -1, 1)$, la transformation $\Phi$ est de longueur 4 donc elle n'est pas à perte maximum. Elle n'est pas paire puisque $\phi$ n'est pas impaire, et elle n'est pas dichotomique puisque $\phi(-1, 1, -1, 1) \neq -\phi(1, -1, 1, -1)$, ce qui est incompatible avec une décomposition de la forme $\phi(x_1, \ldots, x_d) = \gamma x_d \psi(x_1, \ldots, x_{d-1}) \psi(x_2, \ldots, x_d)$.

On ne peut donc pas écrire $\Phi$ comme composée d'une application de longueur 3 et d'une application de longueur 1, puisque les seules applications de longueur 1 sont $\Phi_1$ et $-\Phi_1$. On ne peut pas non plus écrire $\Phi$ comme composée de deux applications de longueur 2, puisqu'une application de longueur 2 perd 0 ou 2 bits. La transformation $\Phi$ est donc indécomposable.

3.5. *Trois exemples de transformations conservant l'information.*   Nous présentons ici trois exemples de transformations homogènes de longueur finie, dont nous montrons qu'elle conservent l'information grâce à la relation d'accordabilité.

PREMIER EXEMPLE.   $H_n = \max(\varepsilon_{n-a_1}, \ldots, \varepsilon_{n-a_m})$ avec $m \geq 2$ et $a_m > \cdots > a_1 \geq 1$.

La transformation est de longueur $d = a_m$ et est associée à l'application $\phi : \{-1, 1\}^d \to \{-1, 1\}$ définie par $\phi(x_0, \ldots, x_{d-1}) = \max(x_{d-a_1}, \ldots, x_{d-a_m})$.

On peut supposer que le PGCD de $a_1, \ldots, a_m$ vaut 1 sans perte de généralité. En effet, si $c$ est un diviseur commun de $a_1, \ldots, a_m$, et si $a_k = cb_k$ pour tout $k \in [1 \ldots m]$, appliquer la transformation

$$(x_n)_{n \in \mathbf{Z}} \mapsto (\max(x_{n-a_1}, \ldots, x_{n-a_m}) x_n)_{n \in \mathbf{Z}}$$

à la suite $(\varepsilon_n)_{n \in \mathbf{Z}}$ revient à appliquer la transformation

$$(x_n)_{n \in \mathbf{Z}} \mapsto (\max(x_{n-b_1}, \ldots, x_{n-b_m}) x_n)_{n \in \mathbf{Z}}$$

aux suites $(\varepsilon_{cn+r})_{n \in \mathbf{Z}}$ pour $r \in [0 \ldots c-1]$.

On utilise alors le fait classique suivant: si une partie de $\mathbf{N}$ contient 0 et est stable pour les translations $n \mapsto n + a_k$, elle contient tous les entiers naturels à partir d'un certain rang $n_0$.

Notons $f_+$ et $f_-$ sont les applications associées à $\phi$. Soit $(x_0, \ldots, x_{d-1}) \in \{-1, 1\}^d$. Définissons la suite $(x_n)_{n \in \mathbf{Z}}$ définie par la relation de récurrence $x_n = \phi(x_{n-d}, \ldots, x_{n-1})$, qui s'écrit aussi $(x_{n-d+1}, \ldots, x_n) = f_+(x_{n-d}, \ldots, x_{n-1})$. Alors l'ensemble des entiers $n \geq d-1$ vérifiant $x_n = 1$ est stable pour les translations $n \mapsto n + a_k$. Si $x_i = 1$ pour un certain $i \in [0 \ldots d-1]$, alors $x_j = 1$ en notant $j = i+d$ si $i < d-1$ et $j = i$ si $i = d-1$. Donc $x_{j+n} = 1$ pour tout $n \geq n_0$, d'où $x_n = 1$ pour tout $n \geq n_0 + 2d - 2$. Par conséquent, $f_+^{n_0 + 2d - 2}$ accorde sur $(1, \ldots, 1)$ tous les $d$-uplets autres que $(-1, \ldots, -1)$.

Comme $f_-$ envoie $(-1, \ldots, -1)$ sur $(-1, \ldots, -1, 1)$ et $(1, \ldots, 1)$ sur $(1, \ldots, 1, -1)$, la composée $f_+^{n_0 + 2d - 2} \circ f_- \circ f_+^{n_0 + 2d - 2}$ accorde sur $(1, \ldots, 1)$ tous les éléments de $E_d$.

DEUXIÈME EXEMPLE. $H_n = \max(\varepsilon_{n-d}, -\varepsilon_{n-c}, \psi(\varepsilon_{n-d}, \ldots, \varepsilon_{n-1}))$ où $d > c \geq 1$ et $\psi$ est une application quelconque de $\{-1, 1\}^d$ dans $\{-1, 1\}$.

La transformation est associée à une application $\phi \colon \{-1, 1\}^d \to \{-1, 1\}$ de la forme $\phi(x_1, \ldots, x_d) = \max(x_1, -x_{e+1}, \psi(x_1, \ldots, x_d))$ avec $e = d - c$. Autrement dit, l'application $\phi$ vérifie la propriété suivante:

Si $x_1 = 1$ ou $x_{e+1} = -1$, alors $\phi(x_1, \ldots, x_d) = 1$.

Montrons que pour $n$ assez grand, la composée $(f_- \circ f_+)^n$ accorde tous les éléments de $E_d$. Pour cela, prenons $(x_1, \ldots, x_d) \in \{-1, 1\}^d$ et définissons une suite $(x_n)_{n \geq 1}$ en posant pour tout $n > d$

$$x_n = \alpha_n \phi(x_{n-d}, \ldots, x_{n-1})$$

avec

$$\alpha_n = \begin{vmatrix} +, & \text{si } (2k-1)d < n \leq 2kd, \\ -, & \text{si } 2kd < n \leq (2k+1)d. \end{vmatrix}$$

Cette relation de récurrence s'écrit aussi $(x_{n-d+1}, \ldots, x_n) = f_{\alpha_n}(x_{n-d}, \ldots, x_{n-1})$. Remarquons que pour tout $n \geq 1$,

$$x_n = 1 \implies x_{n+d} = \alpha_{n+d},$$

$$x_n = -1 \implies x_{n+c} = \alpha_{n+c}.$$

On en déduit immédiatement que les $e$ dernières composantes du $d$-uplet $(x_{2d+1}, \ldots, x_{3d})$ valent $-1$. En effet, pour tout $k \in [1 \ldots e]$:

– si $x_{k+d} = 1$, alors $x_{k+2d} = \alpha_{k+2d} = -1$, d'où $x_{k+2d+c} = \alpha_{k+2d+c} = -1$;
– si $x_{k+d} = -1$, alors $x_{k+d+c} = \alpha_{k+d+c} = 1$, d'où $x_{k+2d+c} = \alpha_{k+2d+c} = -1$.

Le même raisonnement montre plus généralement que pour tout $n \in \mathbf{N}^*$, les $e$ dernières composantes du $d$-uplet $(x_{2nd+1}, \ldots, x_{(2n+1)d})$ valent $-1$.

On utilise ensuite le fait suivant: pour tout $n \in \mathbf{N}$ et pour tout $k \in [e+1 \ldots d]$,

$$x_{2nd+k} = -1 \quad \Longrightarrow \quad x_{(2n+1)d+k-e} = \alpha_{(2n+1)d+k-e} = 1$$
$$\Longrightarrow \quad x_{(2n+2)d+k-e} = \alpha_{(2n+2)d+k-e} = -1.$$

Une récurrence immédiate montre alors que pour tout $n \in \mathbf{N}^*$, les $\min(ne, d)$ dernières composantes du $d$-uplet $(x_{2nd+1}, \ldots, x_{(2n+1)d})$ valent $-1$. Par conséquent, si $n \geq d/e$, la composée $(f_- \circ f_+)^n$ est constante égale à $(-1, \ldots, -1)$.

Troisième exemple.   $H_n = \phi(\varepsilon_{n-d}, \ldots, \varepsilon_{n-1})$ où $\phi \colon \{-1, 1\}^d \to \{-1, 1\}$ est une application différente de l'application $(x_1, \ldots, x_d) \mapsto x_d$ vérifiant $\phi(x_1, \ldots, x_d) \geq x_d$ pour tout $(x_1, \ldots, x_d) \in \{-1, 1\}^d$. Autrement dit:

– pour tout $(x_1, \ldots, x_{d-1}) \in \{-1, 1\}^{d-1}$, $\phi(x_1, \ldots, x_{d-1}, 1) = 1$;
– il existe $(x_1, \ldots, x_{d-1}) \in \{-1, 1\}^{d-1}$ tel que $\phi(x_1, \ldots, x_{d-1}, -1) = 1$.

On voit immédiatement que sous ces hypothèses, l'ensemble des $d$-uplets finissant par 1 est stable pour l'application $f_+$, et que $f_+^d$ envoie tous ces $d$-uplets sur $(1, \ldots, 1)$. Mais $f_+$ envoie au moins un $d$-uplet finissant par $-1$ sur un $d$-uplet finissant par 1, l'application $f_+^{d+1}$ envoie ce $d$-uplet sur $(1, \ldots, 1)$, comme tous les $d$-uplets finissant par 1.

On a donc $2^{d-1} + 1$ éléments simultanément accordables par $f_+^{d+1}$. Comme le nombre maximum d'éléments simultanément accordables divise $2^d$, il vaut donc $2^d$.

3.6. *Questions relatives au décentrage.*  A la fin de la deuxième partie, nous avons montré que si une transformation de la forme $\Phi \colon (\varepsilon_n)_{n \in \mathbf{Z}} \to (H_n \varepsilon_n)_{n \in \mathbf{Z}}$ avec $H_n = \phi(\varepsilon_{n-d}, \ldots, \varepsilon_{n-1})$ vérifie $P[H_n \neq 1] < \frac{1}{2d}$, elle conserve l'information.

De façon plus générale, on s'attend à ce qu'une transformation homogène pour laquelle le processus prévisible $(H_n)_{n \in \mathbf{Z}}$ est très décentré (autrement dit $\mathbf{E}[H_n]$ très proche de 1 ou de $-1$) conserve l'information, puisqu'une telle transformation est "proche" de $\pm \mathrm{Id}$. Nous savons que pour une transformation homogène de longueur $\leq d$, la condition $|\mathbf{E}[H_n]| > 1 - \frac{1}{d}$ assure la conservation de l'information. Une question naturelle est de savoir si on peut remplacer $1 - \frac{1}{d}$ par une borne $r < 1$ indépendante de la longueur.

Nous ne connaissons pas la réponse à cette question, mais nous savons que si une telle borne existe, elle est au moins égale à $\frac{2}{3}$. En effet, il est facile de construire des transformations homogènes perdant de l'information pour lesquelles $|\mathbf{E}[H_n]|$ est arbitrairement proche de $\frac{2}{3}$; on peut même atteindre $|\mathbf{E}[H_n]| = \frac{2}{3}$ avec une transformation homogène de longueur infinie. L'idée pour construire de telles transformations est de regarder parmi les transformations "dichotomiques."

*Un exemple de longueur infinie.* Pour $x = (x_n)_{n \in -\mathbf{N}} \in \{-1; 1\}^{-\mathbf{N}}$, notons $m(x) = \min\{k \in \mathbf{N} : x_{-k} = 1\}$ et $h(x) = (-1)^{m(x)}$, en convenant que $m(x) = 0$ si la suite $(x_n)_{n \in -\mathbf{N}}$ est constante égale à $-1$. Soit $\Psi$ la transformation définie par $\Psi(x) = (h(x_{n-1]})x_n)_{n \in \mathbf{Z}}$ pour tout $x = (x_n)_{n \in \mathbf{Z}} \in \{-1; 1\}^{\mathbf{Z}}$.

La composée $\Phi = \Phi_1 \circ \Psi$ perd au moins un bit d'information et transforme le jeu de pile ou face $(\varepsilon_n)_{n \in \mathbf{Z}}$ en $(H_n \varepsilon_n)_{n \in \mathbf{Z}}$, avec

$$H_n = \varepsilon_{n-1} h(\varepsilon_{n-1]}) h(\varepsilon_{n-2]}).$$

On vérifie immédiatement que sur l'événement $[\varepsilon_{n-1} = -1]$, $m(\varepsilon_{n-2]}) = m(\varepsilon_{n-1]}) - 1$ presque sûrement, d'où $H_n = 1$, tandis que sur l'événement $[\varepsilon_{n-1} = 1]$, $H_n = h(\varepsilon_{n-2]})$, d'où

$$P[H_n \neq 1] = P[\varepsilon_{n-1} = 1; m(\varepsilon_{n-2]}) \text{ impair}] = \tfrac{1}{6}.$$

Par conséquent $\mathbf{E}[H_n] = P[H_n = 1] - P[H_n = -1] = \frac{2}{3}$.

*Un exemple de longueur finie.* Soit $d \geq 2$. Pour obtenir une transformation $\Phi$ de longueur $d$, il suffit de modifier légèrement la définition de $m(x)$ en posant $m(x) = \min\{k \in [0, \ldots, d-2] : x_{-k} = 1\}$, avec la convention $m(x) = d - 1$ si $x_{-k} = -1$ pour tout $k \in [0, \ldots, d-2]$. De cette façon, $m(x)$ est bien une fonction de $(x_{-d+2}, \ldots, x_0)$, donc $\Psi$ est de longueur $d - 1$ et $\Phi$ est de longueur $d$. Avec ce changement, le calcul devient

$$P[H_n \neq 1] = P[\varepsilon_{n-1} = 1; m(\varepsilon_{n-2]}) \text{ impair}]$$
$$+ P[\varepsilon_{n-1} = -1; m(\varepsilon_{n-1]}) = m(\varepsilon_{n-2]}) = d - 1]$$
$$= P[\varepsilon_{n-1} = 1] P[m(\varepsilon_{n-2]}) \text{ impair}] + P[\varepsilon_{n-1} = \cdots = \varepsilon_{n-d} = -1].$$

Or $P[m(\varepsilon_{n-2]}) = k] = \frac{1}{2^{k+1}}$ pour tout $k \in [0 \ldots d-2]$ et $P[m(\varepsilon_{n-2]}) = d - 1] = \frac{1}{2^{d-1}}$.

Si $d = 2q$ avec $q \in \mathbf{N}^*$, on a donc

$$P[H_n \neq 1] = \frac{1}{2} \sum_{l=0}^{q-2} \frac{1}{4^{l+1}} + \frac{1}{2^d} + \frac{1}{2^d} = \frac{1}{2} \frac{4^{q-1} - 1}{3.4^{q-1}} + \frac{2}{2^d} = \frac{1}{6} \frac{2^d + 8}{2^d}.$$

Si $d = 2q + 1$, avec $q \in \mathbf{N}^*$ on a donc

$$P[H_n \neq 1] = \frac{1}{2} \sum_{k=0}^{q-1} \frac{1}{4^{k+1}} + \frac{1}{2^d} = \frac{1}{2} \frac{4^q - 1}{3.4^q} + \frac{1}{2^d} = \frac{1}{6} \frac{2^d + 4}{2^d}.$$

Donc $P[H_n \neq 1]$ est aussi proche de $\frac{1}{6}$ qu'on veut.

*Autres questions.* Une autre question naturelle est de savoir si pour une transformation dichotomique, on a nécessairement $|\mathbf{E}[H_n]| \leq \frac{2}{3}$. Si on reprend le jeu de pile ou face associé au développement dyadique d'un réel défini au Paragraphe 1.2, on peut ramener cette question à un problème d'analyse: notons $\varepsilon_0(t) = -1$ pour $0 \leq t < \frac{1}{2}$ et $\varepsilon_0(t) = 1$ pour $\frac{1}{2} \leq t < 1$. A-t-on pour toute application mesurable 1-périodique $f\colon \mathbf{R} \to \{-1, 1\}$,

$$\left| \int_0^1 f(t) f(2t) \varepsilon_0(t)\, dt \right| \leq \frac{2}{3}?$$

L'inégalité ci-dessus est une égalité lorsque $f$ vaut 1 sur les intervalles $[\frac{1}{2^{2k+1}}, \frac{1}{2^{2k}}[$ pour tout $k \in \mathbf{N}$ et $-1$ ailleurs sur $[0, 1[$. Remarquons que la valeur absolue n'apporte rien de plus à l'inégalité puisque changer $f$ en $\widetilde{f}\colon t \mapsto f(-t)$ transforme l'intégrale en son opposée.

Une dernière question provient de l'observation des résultats obtenus par l'informatique. Pour $d \geq 1$ et $\phi\colon \{-1; 1\}^d \to \{-1; 1\}$, appelons décentrage de $\phi$ la quantité $\text{Card}(\phi^{-1}(1)) - 2^{d-1}$. Lorsque $d \leq 4$, nous avons constaté que toutes les applications $\phi$ pour lesquelles la transformation $\Phi$ associée perd de l'information ont un décentrage pair. Cela reste-il vrai pour toute valeur de $d$?

## RÉFÉRENCES

Institut Fourier
Laboratoire de Mathématiques
UMR5582 (UJF-CNRS)
BP 74
38402 St Martin D'Hères Cedex
France
E-mail: christophe.leuridan@ujf-grenoble.fr
URL: http://www-fourier.ujf-grenoble.fr/~leuridan